\documentclass[10pt]{amsart}

\usepackage{latexsym,enumitem}
\usepackage{amssymb}
\usepackage[cp850]{inputenc}
\usepackage{epsfig}
\usepackage{psfrag}
\usepackage{amsthm}
\usepackage{amscd}
\usepackage{amsmath}
\usepackage{amsfonts}
\usepackage{graphics,caption}
\usepackage[all]{xy}
\usepackage{etoolbox}
\patchcmd{\quote}{\rightmargin}{\leftmargin 2em \rightmargin}{}{}
 \usepackage[sort,numbers]{natbib} 
 \usepackage{mathtools}
\usepackage[normalem]{ulem}

\usepackage[english]{babel}
\usepackage{comment}
\usepackage{color}
\usepackage[colorlinks]{hyperref}

\let\phi\varphi

\DeclareMathOperator{\eqdef}{\coloneqq} 
\newcommand{\f}[2]{\frac{#1}{#2}} 

\let\epsilon\varepsilon
\let\subset\subseteq
\newcommand{\be}{\begin{equation*}}
 \newcommand{\ee}{\end{equation*}}
\newcommand{\bpf}{\begin{dimo}}
\newcommand{\epf}{\end{dimo}}
\newcommand{\bdefi}{\begin{defin}}
\newcommand{\edefi}{\end{defin}}
\newcommand{\bthm}{\begin{thm}}
\newcommand{\ethm}{\end{thm}}
\newcommand{\blem}{\begin{lem}}
\newcommand{\elem}{\end{lem}}
\newcommand{\bcor}{\begin{cor}}
\newcommand{\ecor}{\end{cor}}

\newcommand{\bprop}{\begin{prop}}
\newcommand{\eprop}{\end{prop}}
\newcommand{\bese}{\begin{ese}}
\newcommand{\eese}{\end{ese}}
\newcommand{\brem}{\begin{rem}}
\newcommand{\erem}{\end{rem}}
\newcommand{\bpfc}{\begin{dimoclaim}}
\newcommand{\epfc}{\end{dimoclaim}}

\newcommand{\rar}{\rightarrow} 

\newcommand{\abs}[1]{\left\lvert#1\right\rvert}						
\newcommand{\set}[1]{\left\{#1\right\}}					
\newcommand{\floor}[1]{\left\lfloor#1\right\rfloor}					
	
\DeclareMathOperator{\emp}{\varnothing} 
\DeclareMathOperator{\N}{\mathbb N}			
\DeclareMathOperator{\R}{\mathbb R}			
\DeclareMathOperator{\Z}{\mathbb Z}			

\newcommand{\rank}[1]{\text{rank}(#1)} 

\newcommand{\w}{\omega} 

\newenvironment{quot}
{
	\vspace{-0.2cm}
	\vspace{0.2cm}
}

\theoremstyle{definition}
\newtheorem{d1}{Definition}[section] 

\newenvironment{defin}
{
	\begin{quot}
		\begin{d1}
		}
		{\end{d1}
	\end{quot}

}

\theoremstyle{definition}
\newtheorem{r1}[d1]{Remark}

\newenvironment{rem}
{
	\begin{quot}
		\begin{r1}
		}
		{\end{r1}
	\end{quot}
}

\theoremstyle{definition}
\newtheorem{e1}[d1]{Exercise}

\theoremstyle{definition}
\newtheorem{ese1}[d1]{Example}

\newenvironment{ese}
{
	\begin{quot}
		\begin{ese1}
	}
	{	
		\end{ese1}
	\end{quot}
}

\theoremstyle{definition}

\theoremstyle{definition}
\newtheorem{f2}[d1]{Fact}

\theoremstyle{definition}

\theoremstyle{definition}

\theoremstyle{definition}
\newtheorem{t1}[d1]{Theorem}

\newenvironment{thm}
{
	\begin{quot}
		\begin{t1}}
		{\end{t1}
	\end{quot}
}

\theoremstyle{definition}
\newtheorem*{T1*}{Theorem}

\newenvironment{teor*}
{
	\begin{quot}
		\begin{T1*}}
		{\end{T1*}
	\end{quot}
}

\newenvironment{dimo}
{\begin{proof}[Proof]
	}
	{\end{proof}}

\newenvironment{dimoclaim}{\emph{Proof of Claim:}\;}{\hfill$\square$}

	\theoremstyle{definition}
	\newtheorem{l1}[d1]{Lemma}
	
	\newenvironment{lem}
	{
		\begin{quot}
			\begin{l1}}
			{\end{l1}
		\end{quot}
	}
	\theoremstyle{definition}
	\newtheorem{p1}[d1]{Proposition}
	
	\newenvironment{prop}
	{
		\begin{quot}
			\begin{p1}}
			{\end{p1}
		\end{quot}
	}
	
	\theoremstyle{definition}
	\newtheorem{c1}[d1]{Corollary}
	
	\newenvironment{cor}
	{
		\begin{quot}
			\begin{c1}}
			{\end{c1}
		\end{quot}
	}

		\renewenvironment{abstract}
	{\list{}{\rightmargin\leftmargin}%
		\item[\textbf{Abstract:}]\relax}
	{\endlist}

\newtheorem{innercustomthm}{Theorem}
\newenvironment{customthm}[1]
  {\renewcommand\theinnercustomthm{#1}\innercustomthm}
  {\endinnercustomthm}

 \newtheorem*{Theorem*}{Theorem}
 \newtheorem*{Proposition*}{Proposition}
 \newtheorem*{Lemma*}{Lemma}

\usepackage[margin=1.5in]{geometry}
\usepackage{tikz}

\newcounter{notes}%

\newcommand{\Mod}[1]{\text{Mod}(#1)}

\newcommand{\oMod}[1]{\overline{\text{Mod}(#1)}}

\newcommand{\MCG}[1]{\text{MCG}(#1)}
\newcommand{\MCGc}[1]{\text{MCG}_c(#1)}

\newcommand{\PMCG}[1]{\text{PMCG}(#1)}
\newcommand{\Ends}[1]{\text{Ends}(#1)}

\renewcommand{\P}{\mathcal P}
\renewcommand{\bar}{\overline}
\newcommand{\cord}{\mathcal O_{\N}}
\newcommand{\RE}{\mathcal {RE}}
\newcommand{\Cantor}{\mathfrak C}
\newcommand{\E}{\mathcal E}
\newcommand{\U}{\mathcal U}

\renewcommand{\hat}{\widehat}

\begin{document}
\title[Density results for the modular group of infinite-type surfaces]{Density results for the modular group of infinite-type surfaces}
\author[Y. Chandran and T.Cremaschi]{Y. Chandran and T.Cremaschi}
\thanks{The authors would like to thank Ian Biringer, Jing Tao, and Nick Vlamis for organising the workshop, supported by NSF 1654114, where this project was started. The second author was also supported by MSCA 101107744- DefHyp}

\maketitle
\abstract{In this work we show two results about approximating, with respect to the compact-open topology, mapping classes on surfaces of infinite-type by quasi-conformal maps, in particular we are interested in density results. The first result is that given any infinite-type surface $S$ there exists a hyperbolic structure $X$ on $S$ such that $\PMCG S\subseteq \oMod X$, for $\Mod X$ the set of quasi-conformal homeomorphism on $X$. The second result is that given any surface $S$ with countably many ends then there exists a hyperbolic structure $X$ such that $\MCG S= \oMod X$.}
\section{Introduction}
In recent years the geometry and topology of infinite-type hyperbolic manifolds has seen a surge of interest, for example see \cite{C2018b,PV2018,BFT2023} and references therein. Particularly, surfaces of infinite-type and their mapping class groups have been of particular interest.

This paper explores the basic question of when can we study the mapping classes on infinite-type surfaces by looking at maps that also carry geometric informations. Namely given an infinite-type surface $S$ we want to find a hyperbolic structure $X\in\mathcal T (S)$ for which we can approximate any mapping class $\phi\in \MCG S$ by elements of $\Mod X$. As we see in Section \ref{esesection} given $S$ and $X$ is not hard to build mapping classes $h\in\MCG S$ that are not in $\Mod X$. Thus, the best one can achieve is approximating them in the compact-open topology. 

The basic motivation to study this question is that the interplay between the topology and geometry has been widely successful in the finite-type setting to solve a variety of problems and in particular the Nielsen-Thurston classification of mapping classes. Moreover, we would like to point out that for quasi-conformal mapping classes there is a classification, albeit crude, due to Matsuzaki \cite{Mat03}. Other work on the modular group of big surfaces can be found in \cite{He2023}.

In this paper there are two main results.

\begin{customthm}{A}\label{mainA} Let $S$ be an infinite-type surface. Then, there exists infinitely many components of Teichm\"uller space such that $\PMCG S\subset \oMod X$.
\end{customthm}

and

\begin{customthm}{B}\label{mainB}  Given any infinite-type surface $S$ with countably many ends there exists a hyperbolic structure $X$ such that $\MCG S = \oMod X$.
\end{customthm}

A starting point of both results is a Theorem of Patel-Vlamis, see \cite{PV2018}, that states that $\PMCG S$ is the closure of the compactly supported mapping classes, $\MCGc S$, and handle-shifts. Since for any structure $X$ we have that $\MCGc S\subset \Mod X$ to show Theorem \ref{mainA} it suffices to show that all handle-shifts are in $\oMod X$. This, will be done by introducing the notion of $(A,B,C)$-biflutes and then solving a combinatorial problem, see subsection \ref{ABCbifsec} and Section \ref{PMCGsec}.

To prove Theorem \ref{mainB} the main issue is that for any pair of ends $x$ and $y$ with homeomorphic neighbourhoods we need to realise a quasi-conformal class swapping the two ends. Realising this is the bulk of the work and is done in subsection \ref{Regenussec}. Building these mapping classes is also the main issue to generalising Theorem \ref{mainB} to an arbitrary surface with uncountably many ends. This is why we introduce the concept of surfaces with regular ends, $\RE$ surfaces, in subsection \ref{REsec} and prove the following general result.

\begin{customthm}{C}\label{mainC}   Let $(S,\P)$ have property $\RE$ then there exists infinitely many component of $\mathcal T(S)$ such that for all $X\in \mathcal T(S)$ we have $\MCG S=\oMod X$.
\end{customthm}

Theorem \ref{mainC} will allow us to construct infinite families of examples with uncountably many ends, see subsection \ref{REMCGsec}. Moreover, we also show that the surface, with uncountably many homeomorpshism type of ends, of \cite{MR2022} satisfies the conditions in Theorem \ref{mainC}.

Possible applications could also be in constructing hyperbolic mapping tori of infinite-type via an extension of the double limit Theorem \cite{Mc1996}. This, could tie in with work of the second author \cite{CS2018,C20171,C2018b,C2018c,CDS2021,CP22} and work on mapping tori of end-periodic maps \cite{FKLL2023,FKLL2023b}. 

\textbf{Organization.} In Section \ref{background} we recall some facts on the topology of infinite-type surfaces, Teichm\"uller Theory, and descriptive set theory. In Section \ref{background} we also develop the concept of $(A,B,C)$-biflutes that will be key to showing Theorem \ref{mainA} whose proof is Section \ref{PMCGsec}.

Section \ref{esesection} covers many examples and phenomenas that give the intuition for the definition of $\RE$ surfaces. In Section \ref{Constrsec} we have all our main constructions and we show how to build $\RE$ surfaces. Finally in subsection \ref{REMCGsec} we show that $\RE$ surfaces satisfy $\MCG S=\oMod X$ for some $X\in\mathcal T (S)$. 

\textbf{Acknowledgements.} The first author thanks the hospitality of the University of Luxembourg and of CIRM, where parts of this work have been completed. Both authors thank Ian Biringer, Jing Tao, and Nick Vlamis for organising the workshop in 2022.

\tableofcontents

\section{Background}\label{background}
With $\Cantor$ we will denote the Cantor set which is the unique, up to homeomorphism, compact, metrisable, perfect, and totally disconnected topological space \cite{C1997}.

We now recall various facts that we will need.
\subsection{Infinite-type surfaces} 
In this note a surface means an orientable 2-manifold, possibly with boundary. A surface is of \textit{finite-type} if its fundamental group is finitely generated. Otherwise, we say the surface is of \textit{infinite-type}. A finite-type surface is determined up to homeomorphism by the triple $(g,n,b)$ corresponding to the genus, number of punctures, and boundary components respectively. For infinite-type surfaces there is also a topological classification due to \cite{K1923,R1963}, where now we also need to encode the different ways you can escape to infinity. In the remainder of this subsection, we make this precise and state the classification Theorem. See \cite{AV2020} for an introduction to the theory and basic definitions. We recall the following classification Theorem.

\bthm[Classification Theorem for surfaces \cite{K1923,R1963}]\label{classsurfacethm} Let $S_1$ and $S_2$ be surfaces and let $g_i\in\N\cup\set{\infty}$ be their genus. Then, $S_1$ and $S_2$ are homeomorphic if and only if $g_1=g_2$, $b_1=b_2$ and there is a homeomorphism of the end space respecting genus ends. Moreover, for any pair of nested subsets of $A\subset B\subset \Cantor$ there exists an infinite-type surface $S=S(A,B)$ in which $\Ends{ S}=(A,B)$.\ethm

The natural topology on $\MCG S$ is given by the compact-open topology. In \cite{PV2018} the authors show that:
\bthm[Theorem 4, \cite{PV2018}]\label{shiftgenerating} The set of Dehn twists topologically generate $\PMCG S$ if and only if
$S$ has at most one end accumulated by genus. If $S$ has at least two ends accumulated by genus then the set of Dehn twists together with the set of handle shifts topologically generate $\PMCG S$.
\ethm 

\begin{figure}[htb!]
\def\svgwidth{0.8\textwidth}
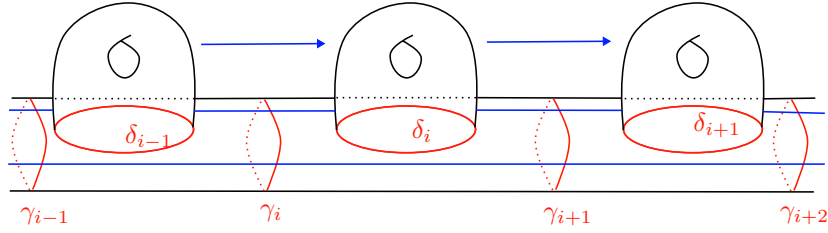\caption{A handle shift on the Ladder surface. The support is in blue and it maps the $\delta_i$ to $\delta_{i+1}$ and $\gamma_i$ to $\gamma_{i+1}$.}\label{handle-shift}
\end{figure}

A handle-shift is the mapping class obtained by shifting countably many handles in a bi-infinite strip, see Figure \ref{handle-shift}.

 We conclude this section with a condition for obtaining complete hyperbolic structures by glueing together pair of pants, see \cite{BS2015}.
\bthm Let $S$ be an infinite-type surface and $\P$ be a pants decomposition. Then, if we have a hyperbolic structure $X\in\mathcal T(S)$ obtained by glueing pants with cuff lengths at most $C<\infty$ and arbitrary twist parameters is complete and of the first kind.\ethm

\subsection{Descriptive Set Theory}
For details see \cite[Chapter 16]{AbSet}. We first recall what countable ordinals are.

 \bdefi We let the set of countable ordinals be denoted by $\cord$. Given an ordinal $\alpha$ we denote by $\mathcal O(\alpha)=\set{\beta\in \cord: \beta<\alpha}$ and by $\overline{\mathcal O(\alpha)}=\set{\beta\in \cord: \beta\leq\alpha}$. We say that a countable ordinal $\w$ is \emph{simple} if it has a unique point of maximal rank, i.e. $\w=\w_0^\zeta+1$ for some ordinal $\zeta$. This is equivalent to saying that its Cantor-Bendixon rank is $(\zeta,1)$.
\edefi 

A key fact is that closed subsets $A \subset \Cantor$ of the Cantor set can be decomposed into a disjoint union $A = B \cup C$ where $B$ is a perfect set (maybe empty), and $B$ is a countable closed set, see \cite[Corollary 1080]{AbSet}.

Roughly, the homeomorphism type of countable closed subsets of $\Cantor$ are determined by the accumulations points and their orders which is made rigorous in the following definition.

\bdefi\label{cantorbendixonrankdefi}
Let $A \subset \Cantor$ be a countable closed  subset of the Cantor set. The \textit{Cantor-Bendixon rank} (or CB rank) of $A$ is the pair $(\nu, n)$ where $\nu$ is the unique (countable) ordinal such that $D^\nu (A) \neq \emptyset$ and $D^{\nu +1}(A) = \emptyset$ and $n = \lvert D^\nu (A) \rvert$. Here $D^\nu$ is $\nu$-th derived set.
\edefi

Here $\nu$ can be any countable ordinal and $n\in \N$ corresponds to number of maximal rank points which is necessarily finite since $\Cantor$ is compact and $D^{\nu +1}(A) = \emptyset$. 

\begin{thm}[Proposition 1086 and Corollary 1095, \cite{AbSet}]
Let $A$ and $B$ be countable closed subsets of the Cantor set, $\Cantor$. Then $A$ is homeomorphic to $B$ if and only if they have the same Cantor-Bendixon rank $(\nu, n)$.
\end{thm}

\brem Each countable ordinal has countably many points but there are uncountably many such ordinals.
\erem 
\bese We start with some easy examples. Any set $E\subset\Cantor$ with finitely many points has finite Cantor-Bendixon rank, in fact it is zero as it has no accumulation points. Their limit is also a countable ordinal, in fact it is the first non-finite countable ordinal $\w_0$. Similarly the set $\w_0\cup\set{\star}$ in which $\star$ is bigger than any element in $\w_0$ is a countable ordinal, the ordinal $\w_0\cup\set\star$ is the first infinite successor ordinal.
\eese

\blem \label{simplecord} Each simple countable ordinal can be constructed recursively as a limit of lower orders one.
\elem 
\bpf See \cite[Theorem 1094]{AbSet}.\epf

We recall the following result from \cite{AbSet}.

\bthm [Theorem 1079 and Corollary 1096, \cite{AbSet}]\label{homeotypecord}Any closed countable subset $C$ of the Cantor set $\Cantor$ is homeomorphic to $\w_0^\zeta n+1$ for $n\in\N$, $\w_0$ the ordinal with order type $\N$, and $\zeta$ a countable ordinal.
\ethm

\subsection{Teichm\"uller theory}
In this section we record some facts from Teichm\"uller Theory that we will use throughout to estimate or bound the quasi-conformal distortions of maps.  

A smooth homeomorphism between Riemann surfaces $f:X\rar Y$ has a derivative which acts as a two dimensional linear map from the tangent space $T_x$ at each point $x \in X$ to the tangent space $T_{f(x)}$. If derivative sends circles to circles then $f$ is \emph{conformal}. Generically the derivative will take a circle in $T_x$ to an ellipse in $T_{f(x)}$. 

\bdefi We say that $f$   is \emph{$K$-quasiconformal} if the eccentricity each of these ellipses is bounded above  by $K$. We will often abbreviate it to qc-map or $K$-qc map.  If additionally $X=Y$ we will say that $f$ is \emph{modular} with respect to $X$. In this context we use the terms modular and quasi-conformal interchangeably. 
\edefi

The following Lemma of Wolpert (see Lemma 12.5 in  \cite{FM2011} ), says that, under $K$-quasiconformal maps, the lengths of hyperbolic geodesics are at most stretched (or compressed) by a factor of $K$. 

\begin{lem}[Wolpert's lemma]\label{wolpert}
Let $f:X \to Y$ be a $K$-quasiconformal homeomorphism between hyperbolic Riemann surfaces. Then for any essential simple closed curve $\alpha$ we have that,
\[\frac{1}{K}\ell_X(\alpha) \leq \ell_Y(f(\alpha))  \leq K \ell_X(\alpha).\]
Where $\ell_X(\alpha)$ denotes the hyperbolic length, in $X$, of the unique geodesic representative isotopic to $\alpha$. 
\end{lem}

The next lemma is a consequence of Theorem 1.1 and Corollary 1.2 in \cite{Bishop} and will be used to normalize the geometry of surfaces for more convenient proofs. We say a hyperbolic surface is of \textit{bounded type} or has \textit{bounded geometry} if for some quantity $M$ it admits a pants decomposition whose lengths are bounded above by $M$ and below by $1/M$.
\blem\label{bishlengths}
Let $X$ be a hyoerbolic surface of bounded type with respect to the pants decomposition $\P$. Then there is a quasiconformal map $f:X \to Y$ where $\ell_Y(f(\alpha)) = 1$ for each $\alpha \in \P$.
\elem

Let $\{\alpha_i\}_{i\in\N}$ be an infinite collection of pairwise disjoint simple closed curves. We consider the mapping class $T = \prod_{i\in\N} T_{\alpha_i}^{n_i}$ where $T_{\alpha_i}$ is the Dehn twists about $\alpha_i$ and $n_i$ is some integer. Matsuzaki gives an estimate on the quasi-conformal dilatation of $T$ acting on a hyperbolic surface $X$ in terms of the lengths of the curves $\alpha_i$ in $X$ (Theorem 1 in \cite{Mat03}). 

\bthm[Theorem 1, \cite{Mat03})]\label{matzmultitwist} 
With the setting as above $T =\prod_{i\in\N} T_{\alpha_i}^{n_i}\in \Mod{X}$. Let $\ell_i$ be the hyperbolic lengths of the geodesic representatives of $\alpha_i$ in $X$. The maximal dilatation $K(T)$ of $T$ acting on $X$ satisfies the following,\small
\[\sup_{i\in\N} \bigg\{ \bigg( \frac{(2 \lvert n_i \rvert -1)_+\ell_i}{\pi} \bigg)^{1/2} +1 \bigg\} \leq K(T) \leq \sup_{i\in\N} \bigg [\bigg\{ \bigg ( \frac{\lvert n_i \rvert \ell_i}{2\theta_i}\bigg )^2 +1 \bigg\}^{1/2} + \frac{\lvert n_i \rvert \ell_i}{2\theta_i} \bigg]^2  \]\normalsize
 where $\theta_i = \pi - 2\arctan\{\sinh(\ell_i/2)\}$ and $(2\lvert n_i \rvert -1)_+ = \max \{(2\lvert n_i \rvert -1), 0\}$.
\ethm
\brem In the setting of Theorem \ref{matzmultitwist} if one has $\f 1 M \leq \ell_i\leq M$ and $\abs{n_i}<C$ then, the induced multi-twist is quasi-conformal. Moreover, one can also boundedly pinch along a multi-curve as long as the pinching amount is uniformly bounded.\erem

\bdefi\label{defiuntwisted} We say that a hyperbolic structure $X$ with bounded geometry is \emph{untwisted} if there exists a bounded pants decomposition $\P$ for $X$ with pants cuff bounded by $\f 1 M\leq \ell_X(c)\leq M$ and such that in Fenchel-Nielsen parameters, given by $\P$, the twists coefficients are uniformly bounded in absolute value.
\edefi

Combining Lemma \ref{bishlengths} and Theorem \ref{matzmultitwist} we get the following corollary which motivates the definition of untwisted.

\bcor \label{coruntwisted}
Let $X$ be a hyperbolic structure untwisted with respect to $\P$. Then $X$ is quasi-conformally equivalent to the structure $Y$ whose Frenchel-Nielsen coordinates with respect to $\P$ have all length parameters set to $1$ and they are glued with no twisting.
\ecor

\brem The condition of being untwisted is not always guaranteed. For example if one takes a flute, with the standard pants decomposition, in which all the non-peripheral curves have length one, then $F$ is clearly of bounded type. Now if one assumes that the twists parameters $t_i$ for the $i$-th loop go to infinity, for example $t_i=i$. Then, there is no untwisted bounded pants decomposition. However, there is always an homeomorphic $F'$ that is untwisted. These two are related by a mapping class element.
\erem

\section{Examples}\label{esesection}
\begin{figure}[htb!]
\def\svgwidth{0.8\textwidth}
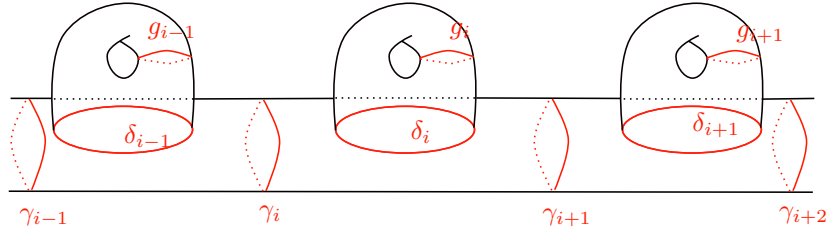\caption{The portion of the ladder $L$ between the $\gamma_{i-1}$ and $\gamma_{i+1}$ separating curves.}\label{figex1}
\end{figure}

If one is already familiar with what potential issues lie in constructing qc-maps of infinite-type surfaces the examples in this section can be skipped. Otherwise, here one can find a collection of examples explaining what can potentially go wrong when trying to realise a mapping class in a modular way. For a precise definition of an $(A,B,C)$-biflute see Definition \ref{ABCbiflutedefin}, for an intuitive definition one can take Figure \ref{figex1} and just have handles every $C$ pants of the flute. Then, the parameters $(A,B)$ depend on how many pants one needs to build the attached $1$-handles.

\bese[A mapping class that is not modular]\label{ese1} Let $L$ be the bi-infinite ladder with the pants decomposition as in Figure \ref{figex1}. We label by $\gamma_i$, $i\in\Z$, the pants curve separating $L$ into two homeomorphic pieces, $\delta_i$ the separating pants curve between $\gamma_i$ and $\gamma_{i+1}$, and $g_i$ the pants curve contained in the compact component obtained by splitting along $d_i$, see Figure \ref{figex1}.

Let $H:L\rar L$ be the full shift mapping each pants curve to the following one of its type, e.g. $\gamma_i\mapsto \gamma_{i+1}$. Since $H$ fixes the two ends we have that $H\in\PMCG L$. 

Let $X_1$ be the hyperbolic structure on $L$ in which all pants curve have length $1$ and there are no twists, this is a complete metric by \cite{BS2015} and is of bounded-type. Then, $H$ maps each pants to an isometric copy itself and so since it is a bijective local isometry (no twists) hence $H$ is an isometry. Thus, $H\in \Mod {X_1}$. 

Let $X_2$ be the structure on $L$ obtained by setting $\ell_{X_2} (\gamma_i)=e^i$ when $i$ is odd, all other curves have length $1$, and no twists. Any properly embedded ray must go out an end and cross infinitely many curves $\gamma_i$ where $i$ is even. The sum of their collars, which give a lower bound for the length such a ray, is infinite so $X_2$ is again complete in the sense of \cite{BS2015}. We now want to show that $H$ is not in $\Mod {X_2}$. 

If it were by Wolpert's Lemma \ref{wolpert} we would have that the dilation $K$ of $H$ satisfies for all $i$:
\[ \f{\ell_{X_2}(H(\gamma_i))}{\ell_{X_2}(\gamma_i)} \leq K,\]
but $H(\gamma_i)=\gamma_{i+1}$ and so when $i$ is even:
\[ \f{\ell_{X_2}(H(\gamma_i))}{\ell_{X_2}(\gamma_i)} =\frac{{e^{i+1}}}{1}=e^{i+1}{\rar}\infty.\]
Moreover, by \cite{PV2018} $H$ is not in $\oMod {X_2}$ which in this case is just the closure of the compactly supported mapping classes.

Also note that $X_2$ does not have any bounded geometry pants decomposition since its length spectrum accumulates to zero.
\eese
\brem In Example \ref{ese1} we have that an element $H\in\MCG L$ such that for a bounded-type metric $X_1$: $H\in\Mod {X_1}$ and for a non bounded-type structure $X_2$ we have that $H\not\in \Mod{X_2}$. Moreover, if one changes the handles with punctures the same construction gives a map that is not in the pure mapping class group and the same arguments can be run.
\erem 

Thus, we saw examples of maps that are not in $\oMod X$ for $X$ that is not of bounded-type but are when $X$ is.

In the following Example \ref{ese2} we show that being bounded-type is not a sufficient condition and we can still have non-modular mapping classes in the mapping class group of a bounded-type surface.

\begin{figure}[htb!]
\def\svgwidth{0.8\textwidth}
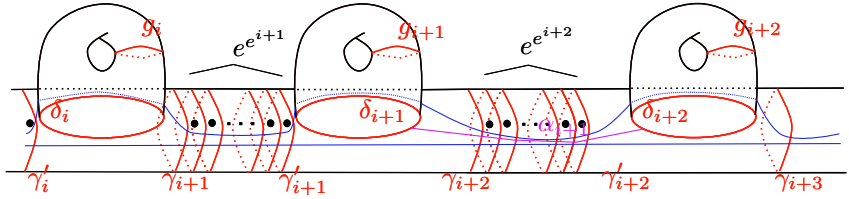
\caption{Topologically $L'$ is the bi-infinite ladder in which both ends are accumulated by punctures. The blue portion is the support of the handle-shift $H$ and one of the loops $\alpha_{i+1}$ which is the outer boundary of a small regular neighbourhood of the complex obtained by adjoining the purple arc, $\delta_{i+1}$, and $\delta_{i+2}$.}\label{figex2.1}
\end{figure}

\bese[Non-modular mapping class in bounded-type]\label{ese2} Consider the surface $L'$ obtained from the bi-infinite ladder $L$ in Example \ref{ese1} in the following way. Let $P_i\subset L$ be the pair of pants bounded by $\gamma_i$, $\gamma_{i+1}$, and $\delta_i$ the loops cutting off $1$-handles as in Figure \ref{figex2.1}.  

We now distance $P_i$ and $P_{i+1}$ by introducing $j_i=e^{e^{\abs {i}}}$ pants $\set{Q_j}_{j=1}^{j_i}$ with a cusp and cuff lengths $1$. Then, the distance between the boundaries of $Q_j$ is some uniform $c>0$. We denote the induced structure by $X_1$ and still by \cite{BS2015} the structure is complete and of bounded-type. Let $H$ be the full handle-shift has before but such that the punctures are not in the support, see Figure \ref{figex2.1}.

Let $\alpha_i$ be the loop going around $\delta_i$ and $\delta_{i+1}$, so that $\alpha_i$ traverses all $Q_j$'s between $P_i$ and $P_{i+1}$ twice. Then, $2c e^{e^i}\leq \ell_{X_1}(\alpha_i)\leq 2c e^{e^i}+C$ for some uniform $C,c>0$ depending on the pants geometry. Still, by Wolpert's Lemma we get that if $H$ was modular its constant $K$ would satisfy:
\[ \frac{e^{e^i(e-1)}} {1 +C/(2c e^{e^i}) } \leq \f{\ell_{X_1}(H(\alpha_i))}{\ell_{X_1}(\alpha_i)}\leq K,\]
which blows up in $i$ giving us a contradiction. 

\begin{figure}[htb!]
\def\svgwidth{0.8\textwidth}
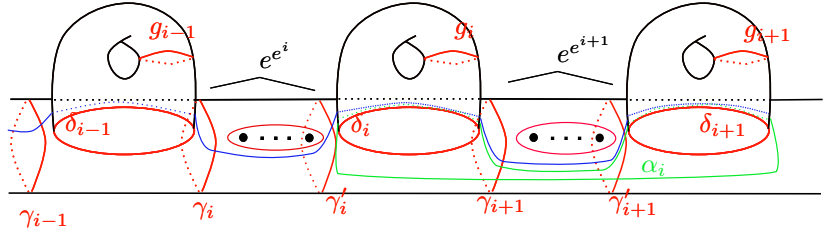
\caption{The changed pants decomposition in which now all punctures are bounded by a pants curve.}\label{figex2.2}
\end{figure}

We now build another structure $X_2$, of bounded-type, on $i$ such that $H$ is modular. To do so we modify the pants decomposition by adding a loop $\beta_i$ that goes around all punctures between $P_i$ and $P_{i+1}$, see Figure \ref{figex2.2}. We call this new pant $Q_i$ and again set all its cuff lengths at $1$.

Now $H$ is QC on $X_2$ and so modular. 
\eese

\brem In Example \ref{ese2} we have an element $H\in\MCG L$ such that it is modular in certain bounded type-structures and not in others. Thus, for example it shows that they are in different Teichm\"uller components. Moreover, as in Example \ref{ese1}, by grabbing, and shifting, a puncture between each $P_i$ and $P_{i+1}$ one can obtain an element with the same properties that is now not a pure mapping class.\erem

\brem
The issue in Example \ref{ese1} and \ref{ese2} is that there are no $(A,B,C)$-biflutes, see Definition \ref{ABCbiflutedefin}. In particular no shift will ever be modular and so in either case we get that the maps are not even in the closure of the modular maps.
\erem 

\bese[pA Example]\label{pAexample} One can also build an example of a mapping class $\psi$ that is never quasi-conformal in any structure on $X$. For simplicity we use the Jacob's ladder $L$ which we now think of as a $\Z$ chain of twice-punctured tori $T_n$, $n\in\Z$.

Let $\psi_n$, $n\in\Z$, be pseudo-Anosov mapping classes on a twice-puncture torus $T$ such that: $\psi_n$ has translation length at least $e^{\abs{n}}$.

Let $\iota_n: T\hookrightarrow S$ mapping $T$ to $T_n$ and let $\psi\eqdef \prod_{n\in\Z} \psi_n\circ \iota_n$. Note that $K(\psi\vert_{T_n}) \leq K(\psi)$ where $K(\cdot)$ is the translation length of a mapping class. Then $\psi$ can never be quasi-conformal since, by (1) and Theorem 2 in \cite{Bers}, $e^{\abs{n}}\leq K(\psi\vert_{T_n})$ which goes to $\infty$ with $n$.

However, note that this is still in the closure of $\oMod X$ since it is in the closure of the compactly supported mapping classes: $\psi=\lim_n \prod_{\abs i<n} \psi_i\circ \iota_i$. Moreover, this construction also gives us that for any $S$ and structure $X$ there are infinitely many homeomorphism that are not quasi-conformal for $X$.
\eese

All examples considered so far had genus and the maps are in the pure mapping class groups. We now show how we have the issues on planar surfaces and with maps that do not fix the end-space.

Let $S_i$ be the planar surface with one boundary component and end space given by a point $x_\infty$ of Cantor-Bendixon rank $i$. So that $S_1$ is just a flute surface with a boundary component.

\bese[Non-modular with only planar ends and bounded-type]\label{ese3}Consider a surface $\Sigma_0$ with countably many boundary components indexed by $\Z$ accumulating on the two ends and equip $\Sigma_0$ with a pants decomposition $\P$ such that each $P\in\P$ contains a unique element of $\partial \Sigma$.

\begin{figure}[htb!]
\def\svgwidth{0.7\textwidth}
\begingroup%
  \makeatletter%
  \providecommand\color[2][]{%
    \errmessage{(Inkscape) Color is used for the text in Inkscape, but the package 'color.sty' is not loaded}%
    \renewcommand\color[2][]{}%
  }%
  \providecommand\transparent[1]{%
    \errmessage{(Inkscape) Transparency is used (non-zero) for the text in Inkscape, but the package 'transparent.sty' is not loaded}%
    \renewcommand\transparent[1]{}%
  }%
  \providecommand\rotatebox[2]{#2}%
  \newcommand*\fsize{\dimexpr\f@size pt\relax}%
  \newcommand*\lineheight[1]{\fontsize{\fsize}{#1\fsize}\selectfont}%
  \ifx\svgwidth\undefined%
    \setlength{\unitlength}{429.23804576bp}%
    \ifx\svgscale\undefined%
      \relax%
    \else%
      \setlength{\unitlength}{\unitlength * \real{\svgscale}}%
    \fi%
  \else%
    \setlength{\unitlength}{\svgwidth}%
  \fi%
  \global\let\svgwidth\undefined%
  \global\let\svgscale\undefined%
  \makeatother%
  \begin{picture}(1,0.25746572)%
    \lineheight{1}%
    \setlength\tabcolsep{0pt}%
    \put(0,0){\includegraphics[width=\unitlength,page=1]{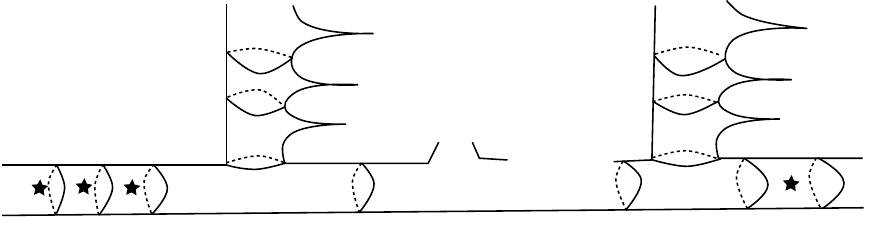}}%
    \put(0.17118859,0.0029928){\color[rgb]{0,0,0}\makebox(0,0)[lt]{\lineheight{1.25}\smash{\begin{tabular}[t]{l}$a_0$\end{tabular}}}}%
    \put(0.83340972,0.00490258){\color[rgb]{0,0,0}\makebox(0,0)[lt]{\lineheight{1.25}\smash{\begin{tabular}[t]{l}$a_1$\end{tabular}}}}%
    \put(0,0){\includegraphics[width=\unitlength,page=2]{images/fig5.pdf}}%
    \put(0.95816281,0.04663109){\color[rgb]{0.97254902,0.01176471,0.1372549}\makebox(0,0)[lt]{\lineheight{1.25}\smash{\begin{tabular}[t]{l}$H$\end{tabular}}}}%
  \end{picture}%
\endgroup%

\caption{The surface $S$ in which we glued $S_1$'s surfaces between the $a_i$, $i\geq 0$, boundary components of the bi-infinite flute $\Sigma_0$. The support of the shift $H$ in the example is red. The stars represent punctures}\label{figex3} 
\end{figure}

To each $a_i\eqdef \floor{ e^{e^i}}-2$, $i\geq 0$, we glue an isolated puncture and fill the boundary components between $a_i$ and $a_{i+1}$ with copies of $S_{i+1}$. Fill all the boundary components of $\Z_{<0}$ with punctures as well, call the resulting surface $S$, see Figure \ref{figex3}. 

Let $H$ be the shift whose support contains all the punctures and shifts the $j$th puncture to the $j+1$, if $j\leq -1$, to $a_0$ if $j=-1$ or $a_j$ to $a_{j+1}$ if $j\geq0$, see Figure \ref{figex3}. Then, as in Example \ref{ese2} we get that $H$ is not quasi-conformal in the complete structure with isometric pants glued with no twists. As Example \ref{ese2} by changing the pants decomposition on $\Sigma_0$ and using isometric pants we get that $H$ is modular.\eese

We conclude by showing that even if we have quasi-conformal shifts between every pair of genus ends we can still have non-conformal mapping classes.

\begin{figure}[htb!]
\def\svgwidth{0.7\textwidth}
\begingroup%
  \makeatletter%
  \providecommand\color[2][]{%
    \errmessage{(Inkscape) Color is used for the text in Inkscape, but the package 'color.sty' is not loaded}%
    \renewcommand\color[2][]{}%
  }%
  \providecommand\transparent[1]{%
    \errmessage{(Inkscape) Transparency is used (non-zero) for the text in Inkscape, but the package 'transparent.sty' is not loaded}%
    \renewcommand\transparent[1]{}%
  }%
  \providecommand\rotatebox[2]{#2}%
  \newcommand*\fsize{\dimexpr\f@size pt\relax}%
  \newcommand*\lineheight[1]{\fontsize{\fsize}{#1\fsize}\selectfont}%
  \ifx\svgwidth\undefined%
    \setlength{\unitlength}{314.8636441bp}%
    \ifx\svgscale\undefined%
      \relax%
    \else%
      \setlength{\unitlength}{\unitlength * \real{\svgscale}}%
    \fi%
  \else%
    \setlength{\unitlength}{\svgwidth}%
  \fi%
  \global\let\svgwidth\undefined%
  \global\let\svgscale\undefined%
  \makeatother%
  \begin{picture}(1,0.49885927)%
    \lineheight{1}%
    \setlength\tabcolsep{0pt}%
    \put(0,0){\includegraphics[width=\unitlength,page=1]{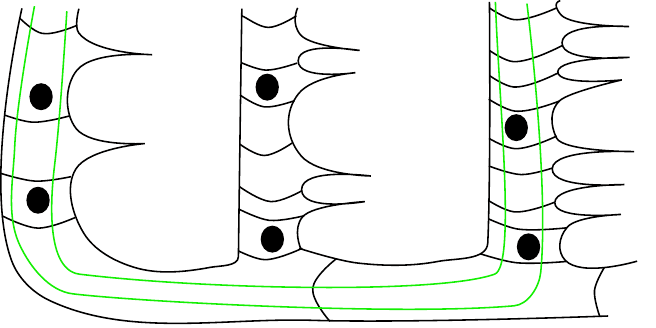}}%
    \put(0.35564564,0.006664){\color[rgb]{0.16470588,0.93333333,0.01176471}\makebox(0,0)[lt]{\lineheight{1.25}\smash{\begin{tabular}[t]{l}$h_{1,3}$\end{tabular}}}}%
    \put(0,0){\includegraphics[width=\unitlength,page=2]{images/fig6.pdf}}%
    \put(0.59517023,0.08346804){\color[rgb]{0.04313725,0.17647059,0.97254902}\makebox(0,0)[lt]{\lineheight{1.25}\smash{\begin{tabular}[t]{l}$h_{2,j}$\end{tabular}}}}%
  \end{picture}%
\endgroup%

\caption{The structure of $S$ with pants decomposition in red and with the support of two shifts $h_{1,3}$ and $h_{2,j}$ drawn in green. The black disks are standard one-handles.}\label{figex5}
\end{figure}

\bese[Non-modular but $(A,B,C)$-biflutes between every pair of ends.]\label{example4} Let $S$ be a surface with end space given by a point of Cantor-Bendixon rank two accumulated by rank one genus ends $E_i$ and punctures. The rank genus ends are as in Figure \ref{figex2.1} with the standard handles separated by $i$ punctures. See Figure \ref{figex5}.

Let $\mathcal P$ be the pants decomposition of Figure \ref{figex5}. Build a complete hyperbolic structure $X$ by glueing pants of type $P\cong P(\ell,\ell,\ell)$ and $Q\cong P(\ell,\ell,0)$ with no twist. Then, the full shift $h_{i,j}$ between the end $E_i$ and end $E_j$ is a modular map of constant $K(i,j)$. Moreover, if either one of $i$ or $j$ goes to infinity then $K(i,j)\rar \infty$ as well. Consider, the pure mapping class:
\[h\eqdef \prod_{i\in\N} h_{i,i+2},\]
since $K(h_{i,i+1})\geq C i$ for some $C>0$ then $h$ is not modular. The map $h$, is well defined since each handle gets moved at most once and so the composition converges in the compact open-topology. However, since each $h_{i,i+2}$ is modular it follows by construction that $h\in\oMod X$.\eese

\section{Modular maps and $(A,B,C)$-biflutes}\label{ABCbifsec}

In this section we build some of the tools we will need for the main result. Namely we introduce $(A,B,C)$-biflutes, and show that, under certain conditions, they generate quasi-conformal maps. Specifically we will show that handle-shifts along bounded geometry $(A,B,C)$-biflutes are quasi-conformal, see Lemma \ref{boundedshiftismodular}.

We recall that for any surface $S$ a pants decomposition $\P$ can be represented by a trivalent graph, where each pair of pants is a vertex, and edges correspond to curves in the decomposition.  By an abuse of notation we will often confuse a pants decomposition and its corresponding graph.

We now define various graphs, and the corresponding pants decomposition, that we'll use throughout this note.

\bdefi \label{graphpantsdecdefi}
A \emph{$\Z$-graph} is the trivalent half-edge graph in which vertices are indexed by $\Z$ and each vertex $v_i$ is glued to $v_{i-1}$ and $v_{i+1}$ and has a free half-edge, the subgraph with vertices $v_i$, $n\leq i \leq m$, is an $[n,m]$-graph. We will use $\P_{\Z}$ for the pants decomposition and $\P_{[n,m]}$ for the corresponding pants decompositions.  Then, $\P_{\Z}$ is the standard pants decomposition on a biflute, and $\P_{[n,m]}$ is a sphere with $m-n+3$ boundaries.

\begin{figure}[htb!]
\def\svgwidth{0.7\textwidth}
\begingroup%
  \makeatletter%
  \providecommand\color[2][]{%
    \errmessage{(Inkscape) Color is used for the text in Inkscape, but the package 'color.sty' is not loaded}%
    \renewcommand\color[2][]{}%
  }%
  \providecommand\transparent[1]{%
    \errmessage{(Inkscape) Transparency is used (non-zero) for the text in Inkscape, but the package 'transparent.sty' is not loaded}%
    \renewcommand\transparent[1]{}%
  }%
  \providecommand\rotatebox[2]{#2}%
  \newcommand*\fsize{\dimexpr\f@size pt\relax}%
  \newcommand*\lineheight[1]{\fontsize{\fsize}{#1\fsize}\selectfont}%
  \ifx\svgwidth\undefined%
    \setlength{\unitlength}{323.79834032bp}%
    \ifx\svgscale\undefined%
      \relax%
    \else%
      \setlength{\unitlength}{\unitlength * \real{\svgscale}}%
    \fi%
  \else%
    \setlength{\unitlength}{\svgwidth}%
  \fi%
  \global\let\svgwidth\undefined%
  \global\let\svgscale\undefined%
  \makeatother%
  \begin{picture}(1,0.17939003)%
    \lineheight{1}%
    \setlength\tabcolsep{0pt}%
    \put(0,0){\includegraphics[width=\unitlength,page=1]{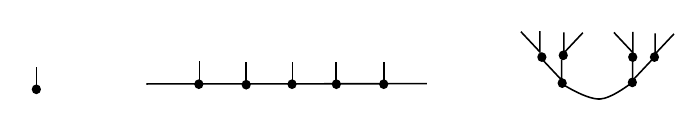}}%
    \put(0.88360625,0.00605939){\color[rgb]{0,0,0}\makebox(0,0)[lt]{\lineheight{1.25}\smash{\begin{tabular}[t]{l}$T$\end{tabular}}}}%
    \put(0,0){\includegraphics[width=\unitlength,page=2]{images/fig1.pdf}}%
  \end{picture}%
\endgroup%

\caption{A loop graph of parameters $(2,3)$, a $[2,6]\subset \Z$ graph, and the regular trivalent graph. The pants decomposition can be obtained by taking the boundary of a  regular neighbourhood in $\R^3$ of the graph and then putting a pair of pants at every vertex. }\label{graphpic}
\end{figure}

A \emph{loop graph} of parameters $(a,b)$, $a\geq 0,b\geq 0$, is the graph with $a+b$ vertices indexed by $1\leq j\leq a+b$. The first $1\leq i \leq a$ vertices are glued $v_i$ to $v_{i+1}$ with $v_1$ and $v_a$ having two free edges. Then, we make another string with the vertices $a+1\leq i\leq a+b$ vertices and we glue $v_{a+1}$ and $v_{a+b}$ to $v_a$ along the previously free edges. If $b=0$ we just glue together the two free edges of $v_a$. If we need to specify the pants decomposition we will use $\P_{\ell(a,b)}$. For examples see Figure \ref{graphpic}.
\edefi 

We will also refer to the pants decomposition on the flute surface induced by $\N$ as $\P_F$ and $\P_T$ will be the pants decomposition for the regular trivalent graph. See Figure \ref{graphpic} for examples of these graphs and pants decompositions.

\brem[Standard Handle]\label{standardhandle} Note that a $(0,1)$-loop graph corresponds to a punctured torus with meridional pants decomposition. This will also be referred to as a \emph{standard handle}.
\begin{figure}[htb!]
\def\svgwidth{0.5\textwidth}
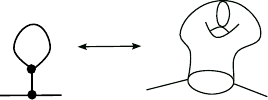
\caption{A $(0,1)$-loop graph attached to a pants and the induced standard one-handle.}\label{stanhandlepic}
\end{figure}
\erem 

There are potentially many choices for loops in the $1$-handle in Figure \ref{stanhandlepic}, namely any primitive $(p,q)$-curve on the punctured torus, however since they all differ by a compactly supported mapping class we will always assume the meridional choice as in Figure \ref{stanhandlepic}, which by Remark \ref{standardhandle} is a standard handle.

\bdefi\label{ABCbiflutedefin} An \emph{$(A,B,C)$-biflute} is an infinite-type surface with at most countably many boundary components, possibly infinite genus and two ends. The surface is defined by the following construction via pants decomposition. Start with pants decomposition given by the graph $\Z$ such that:
\begin{itemize} 
\item every $ C$ places we add a loop graph $L_i$ to the free boundary component;
\item each loop graph $L_i$ is of type $(A,B)$ and if $(0,0)$ no loop graph.
\end{itemize}
We say that a biflute is a \emph{generalised} \emph{$(A,B,C)$-biflute} if the $(A,B,C)$ are not constant but uniformly bounded. That is, we add a loop graph $L_i$ every $C_i\leq C$ places and $L_i$ is an $(A_i,B_i)$ loop graph with $A_i\leq A$,  $B_i\leq B$ and we have that either $A_i + B_i \geq 1$ for all but finitely many $i$, or $A_i + B_i = 0$ for all but finitely many $i$.

Similarly, we define $(A,B,C)$-flutes and generalised $(A,B,C)$-flutes if the base graph is indexed by $\N$.
\edefi

Note that two $(A,B,C)$ and $(A',B',C')$ biflutes (or flutes) are homeomorphic if $A+B = A'+B' = 0$ or $\min\{A+B, A' + B'\} > 0$, however they are combinatorially different and in particular once we put a metric on them they might not even be quasi-conformally equivalent, see the Examples from Section \ref{esesection}.

\bese If one consider the ladder surface $L$ then $L$ can be given the structure of a $(0,1,1)$-biflute, that is we attach a standard one-handle to every free edge of $\P_{\Z}$. See Figure \ref{figex1}.
\eese 

\begin{figure}[htb!]
\def\svgwidth{0.8\textwidth}
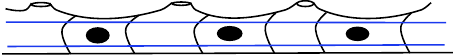
\caption{A shift along a $(0,1,2)$-biflute the black disks are standard one-handles.}\label{shiftsupport}

\end{figure}

\bdefi
Given $F$ an $(A,B,C)$-biflute we define a \emph{shift along $F$} as the full handle-shift with the support as in Figure \ref{shiftsupport}.

\edefi

A key, result for us will be the following observation.

\blem\label{pantsQC} Let $f:(S,\P)\rar (S,\P)$ be a homeomorphism mappings pants to pants. Then, for any structure $X$ in $\mathcal T(S)$ for which both $X$ and $f\cdot X$ have bounded geometry and are untwisted with respect to the pants $\P$ then $f$ is in $\Mod X$.
\elem
\bpf Since both $X$ and $f\cdot X$ are of bounded type and untwisted with respect to $\P$, they are both quasi-conformally equivalent to a structure $Y$ whose Fenchel-Nielsen coordinates with respect to $\P$ have all length parameters set to $1$ and all twist parameters set to $0$ by Corollary \ref{coruntwisted}. Since both $X$ and $f\cdot X$ are quasi-conformally equivalent to $Y$, they are also quasi-conformally equivalent to each other.
\epf

By Lemma \ref{pantsQC} we obtain the following corollary.

\bcor\label{boundedshiftismodular} If $X$ has bounded geometry, with pants decomposition $\P$, then any shift $h$ along an untwisted generalised $(A,B,C)$-biflute is modular.\ecor

\brem \label{allshiftsaremodularifoneis}This will be enough for us since if $h$ is a quasi-conformal shift between the ends $x$ and $y$ and $h'$ is any shift between $x$ and $y$ then we have that $h'\in\overline{\langle h,\MCGc S\rangle}$. The converse is not true. For example take the ladder surface with the $(0,1,1)$ pants decomposition and the full shift $h$ and denote by $\gamma_i$, $i\in\Z$, the separating pants curves leaving the two ends. Let the metric $X$ be such that $\gamma_i=1$ for $i>0$ and $\gamma_i=e^i$ for $i\leq 0$ and all other pants curve have length $1$. Then, $h$ is quasi-conformal, but there is no pants decomposition of bounded geometry realizing the $(A,B,C)$-biflute. \erem

Lemma \ref{boundedshiftismodular} will allows us to treat the problem of building quasi-conformal maps on a surface $S$ combinatorially.

\begin{figure}[htb!]
\def\svgwidth{0.8\textwidth}
\begingroup%
  \makeatletter%
  \providecommand\color[2][]{%
    \errmessage{(Inkscape) Color is used for the text in Inkscape, but the package 'color.sty' is not loaded}%
    \renewcommand\color[2][]{}%
  }%
  \providecommand\transparent[1]{%
    \errmessage{(Inkscape) Transparency is used (non-zero) for the text in Inkscape, but the package 'transparent.sty' is not loaded}%
    \renewcommand\transparent[1]{}%
  }%
  \providecommand\rotatebox[2]{#2}%
  \newcommand*\fsize{\dimexpr\f@size pt\relax}%
  \newcommand*\lineheight[1]{\fontsize{\fsize}{#1\fsize}\selectfont}%
  \ifx\svgwidth\undefined%
    \setlength{\unitlength}{216.38550766bp}%
    \ifx\svgscale\undefined%
      \relax%
    \else%
      \setlength{\unitlength}{\unitlength * \real{\svgscale}}%
    \fi%
  \else%
    \setlength{\unitlength}{\svgwidth}%
  \fi%
  \global\let\svgwidth\undefined%
  \global\let\svgscale\undefined%
  \makeatother%
  \begin{picture}(1,0.39299938)%
    \lineheight{1}%
    \setlength\tabcolsep{0pt}%
    \put(0,0){\includegraphics[width=\unitlength,page=1]{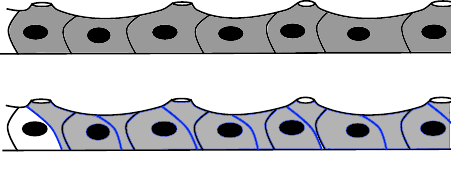}}%
    \put(0.34181941,0.22474629){\color[rgb]{0.6,0.6,0.6}\makebox(0,0)[lt]{\lineheight{1.25}\smash{\begin{tabular}[t]{l}$F$\end{tabular}}}}%
    \put(0.34181941,0.02284263){\color[rgb]{0.6,0.6,0.6}\makebox(0,0)[lt]{\lineheight{1.25}\smash{\begin{tabular}[t]{l}$s(F)$\end{tabular}}}}%
  \end{picture}%
\endgroup%

\caption{A flute shear $s:U\hookrightarrow U$ for a $(0,1,2)$-flute $F$, in grey in the picture. The image of the non-handle pants curve are in blue and the black disks are standard one-handles.} \label{fluteshearfig}
\end{figure}

In the latter part of this work we will also care about shifts along flutes and we define them as generalised biflutes, see Figure \ref{fluteshearfig}.

\bdefi\label{fluteshear} Let $(F,\P_F)\subset (S,\P)$ be the embedding of an $(A,B,C)$-flute, possibly a generalised flute. Let $U\subset S$ be the non-compact surface in $\overline{S\setminus \partial F}$ containing $F$. We define the \emph{shear along $F$} to be the embedding $s:U\hookrightarrow U$ obtained by moving each pants in $F$ to the next one. See Figure \ref{fluteshearfig}.
\edefi 

As for shifts along biflutes we get:
\blem\label{boundedfluteshiftismodular} If $X$ has bounded geometry, with pants decomposition $\P$, and is untwisted then any flute-shear $s$ along a generalised $(A,B,C)$-flute $F$ is modular. \elem
\bpf Without loss of generality we can assume that each pant in $F$ is isometric to a pant with cuff lengths one. Since $s$ is the identity outside of $F$ we only need to check that $s$ is quasi-conformal on $F$. Recall that the shear $s:F\rar F$ is defined as in Figure \ref{fluteshearfig}.

On the level of the geometry a single block, in the case of a $(0,1,2)$-flute shear, gets stretched as in the following Figure \ref{blockshear}.
\begin{figure}[htb!]
\def\svgwidth{0.8\textwidth}
\begingroup%
  \makeatletter%
  \providecommand\color[2][]{%
    \errmessage{(Inkscape) Color is used for the text in Inkscape, but the package 'color.sty' is not loaded}%
    \renewcommand\color[2][]{}%
  }%
  \providecommand\transparent[1]{%
    \errmessage{(Inkscape) Transparency is used (non-zero) for the text in Inkscape, but the package 'transparent.sty' is not loaded}%
    \renewcommand\transparent[1]{}%
  }%
  \providecommand\rotatebox[2]{#2}%
  \newcommand*\fsize{\dimexpr\f@size pt\relax}%
  \newcommand*\lineheight[1]{\fontsize{\fsize}{#1\fsize}\selectfont}%
  \ifx\svgwidth\undefined%
    \setlength{\unitlength}{98.5140042bp}%
    \ifx\svgscale\undefined%
      \relax%
    \else%
      \setlength{\unitlength}{\unitlength * \real{\svgscale}}%
    \fi%
  \else%
    \setlength{\unitlength}{\svgwidth}%
  \fi%
  \global\let\svgwidth\undefined%
  \global\let\svgscale\undefined%
  \makeatother%
  \begin{picture}(1,0.15779683)%
    \lineheight{1}%
    \setlength\tabcolsep{0pt}%
    \put(0,0){\includegraphics[width=\unitlength,page=1]{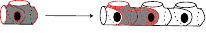}}%
    \put(0.08557411,0.01510087){\color[rgb]{0.05882353,0,0}\makebox(0,0)[lt]{\lineheight{1.25}\smash{\begin{tabular}[t]{l}$\gamma$\end{tabular}}}}%
    \put(0.65354026,0.01510087){\color[rgb]{0.05882353,0,0}\makebox(0,0)[lt]{\lineheight{1.25}\smash{\begin{tabular}[t]{l}$s(\gamma)$\end{tabular}}}}%
  \end{picture}%
\endgroup%

\caption{For the (0,1,2)-flute shear $s:F\hookrightarrow F$. Each pants curve on the left, in red, is mapped to a pants curve on the right except for $\gamma$.} \label{blockshear}
\end{figure}
Since the geometry is the same everywhere each such block gets distorted a uniform amount and so the whole map is quasi-conformal. 
\epf

We conclude with the key Lemma about surfaces with biflutes. This will allow us to show that $\PMCG S\subset \oMod X$ for $X$ a finite-type surface over $S$. 

\bdefi\label{abccond} An infinite-type surface  $S$ with a pants decomposition $\P$ satisfies the \emph{biflute} condition if for every pair of genus ends $e_1$ and $e_2$ there is some $(A,B,C)$-biflute.  Moreover, if $S$ has boundary we also require that any genus end has an $(A,B,C)$-flute to every boundary component.
\edefi 
By Lemma \ref{boundedshiftismodular} and Remark \ref{allshiftsaremodularifoneis} we obtain that for all surfaces with $(A,B,C)$-biflutes modular shifts contain all shifts in $\Mod X$.

\bcor\label{bifluteconditionforQC} If $X$ has bounded geometry with untwisted pants decomposition $\P$ and satisfies the biflute condition then there is a modular shift between every pair of non-planar ends and a modular shear between every boundary component and every genus end. Moreover, any shift is in $\oMod X$.\ecor
\bpf The first part of the statement follows from Lemma \ref{boundedshiftismodular} so we only need to show that any handle-shift $h$ is in $\oMod X$. Let $h$ be a shift between the two genus ends $x$ and $y$ and let $h'$ be the modular handle shift between $x$ and $y$ coming from the biflute condition and Lemma \ref{boundedshiftismodular}. Then, since $h$ can be obtained as the limit of $f_n\circ h'$ with $f_n$ compactly supported the result follows, see Remark \ref{allshiftsaremodularifoneis}. \epf

\section{For any countable surface $S$ there exists a structure $X$ such that the Modular group contains $\PMCG S$}\label{PMCGsec}

The aim of this section is to show that given any infinite-type surface $S$ we have $\PMCG S\subset\oMod X$ for a suitable choice of structure $X$ on $S$. By the result of Patel-Vlamis \cite{PV2018} and the previous section we can solve this by showing that any $S$ admits a pants decomposition with generalised $(A,B,C)$-biflutes, this is done in Proposition \ref{pantsforpure}.

To us the following surfaces will be key in the construction:

\bdefi\label{genflute} We say that a surface $S=S_\eta$, $\eta\in\cord$, is a \emph{$\eta$-sphere} if it is planar and its end space has a unique point of order $\eta$, that is $\Ends S\cong \w_0^\eta+1$.
\edefi 
\bese[$\eta$-spheres]\label{genfluteese}
Then, $S_1$ is the flute surface, and $S_0$ is a sphere with a puncture. The surfaces $S_i$ in Example \ref{ese3} are $i$-spheres. By the Classification Theorem \ref{classsurfacethm} and the ordinal structure, see Lemma \ref{simplecord}, an $\eta$-sphere can be built inductively by taking a boundary flute and glueing to each boundary component a surface $S_{\eta_i}$ with $\eta_i\rar \eta$.
\eese

\bdefi Given a topological space $T$ with a marked point $p\in T$ and a connected separating set $\gamma$ we call $B_\gamma$ the \emph{branch at $\gamma$} to be the component of $T\setminus \gamma$ not containing $p$.
\edefi 

We now introduce the concept of level structure, this is inspired by the natural level structure of a regular tree.

\bdefi
Let $S$ be a surface with a pants decomposition $\P$. A \emph{level structure} $\Lambda$ on $S$ is a sub-collection, not necessarily proper, of loops of $\P$ such that:
\begin{enumerate}
\item each $\gamma\in\Lambda$ is separating;
\item $\Lambda=\coprod_{n=0}^\infty \Lambda_n$ such that splitting along $\Lambda_n$ gives a compact sub-surface $X_n$, and all $\Lambda_j$, $j\geq n+1$, are in $X_n^C$;
\item if $S$ has finitely many boundary components then $\Lambda_0$ is comprised only of the boundary components of $S$;
\item components of $\overline{X_{n+1}\setminus X_n}$ are either pair of pants, twice or thrice punctured tori, or annuli exhausting an isolated end of $S$.
\end{enumerate} 
\edefi

\bese[Level Structures] \label{levelstructure} Consider the Cantor tree with the pants decomposition $\P_T$ given by the regular trivalent graph. Then, $\Lambda_T=\P_T$ is a level structure in which $\Lambda_n$ is the collection of $2^n$ loops at level $n$.

Similarly, consider the $i$-spheres $S_i$'s with one boundary component. Let $F$ be the flute with the $\N$ pants decomposition $\P_{\N}$ in which $\Lambda_0$ is the left-most boundary and the boundaries attached to level $i$ have level $i+1$. Thus, going from left to right in $F$ the boundaries component are at level $i\geq 0$. This gives $S_1$. Then, we continue recursively. The surface $S$ obtained by glueing a $S_i$ surface to the $i$-th boundary component of $F$ also has a level structure $\Lambda(S)$ obtained by taking, set-wise, $\Lambda(F)\cup_i\Lambda(S_i)$. In which the original levels of $\Lambda(S_i)$ are shifted by $i$. That is, 
\[ \Lambda_k(S)=\Lambda_k(F)\cup \cup_{i<k}\Lambda_{k-i}(S_i).\]
This construction and the fact that $S_i$ is obtained by gluing $S_{i-1}$ along a flute $F$ gives a level structure on any $S_i$. The same works for any $\eta$-sphere. \eese 

Using the level structure of the regular tree we can now prove the main topological result for showing that $\PMCG S\subset \oMod X$ for a carefully chosen $X$. Recall that a $(0,1,2)$-biflute is as in Figure \ref{shiftsupport}.

\bprop\label{pantsforpure} Let $S$ be an infinite-type surface with no non-compact boundary components, then $S$ admits a pants decomposition $\P$ such that each pair of genus ends has a $(0,1,2)$-biflute.
\eprop 
\bpf We first delete all boundary components to obtain a new surface, which we still denote by $S$, in which some isolated planar ends correspond to boundaries, we call this collection $\mathfrak b\subset \Ends S$. 

The end space of $S$ is $\Ends S=(\E,\E_G)$. Let $\iota: \Ends S\hookrightarrow \Cantor$ be an embedding into the Cantor set thought as the end space the binary Cantor Tree $S_T$, which we assume to be rooted at $p$. Let $\P_T$ be the standard pants decomposition on the Cantor Tree as the regular trivalent graph. Recall that for $\P_T$ we have a level structure $\Lambda=\cup_{i\in\N} \Lambda_i$ in which at level $i$ we have $2^i$ loops.

We then alter $(S_T,\P_T)$ in the following way. Each $\gamma \in \pi_0(\P_T)$ gives us a branch $B_\gamma$ of the tree. If $\iota(\Ends S)\cap \Ends{ B_\gamma}=\emp$ we delete the whole branch $B_\gamma$, replace it with a disk and remove $\gamma$ from $\P_T$. 

While doing this procedure the following can happen. After replacing $\gamma\in\Lambda_i$ with a disk some $\beta\in\Lambda_i$ can become isotopic to $\alpha\in\Lambda_{i-1}$. If this is the case, and $B_\beta$ has more than one point of $\iota(\E)$ we delete $\beta$ from $\P_T$ and decrease the level by one of all loops $\delta\in B_\beta\cap \P_T$. If $\iota(\E)\cap \Ends{B_\beta}$ has only one point we delete the whole branch and glue in $\mathbb S^1\times [0,\infty)$ in which at the integers point we mark a loop $\gamma_j=\mathbb S^1\times \set j$ which we start at level $i+j$.

By doing this process inductively level by level for all pants curve we obtain a surface $(S_0,\P_0)$ such that:
\begin{enumerate}
\item $\Ends{S_0}=\E$;
\item $\P_0$ is a pants decomposition on $S_0$ in which each curve is separating and if two loops are isotopic then they are in an isolated end of $S_0$;
\item $\P_0$ inherits a level structure from $\P_T$ however each level has most $\ell_i$ separating curves, $\ell_i\leq 2^i$, no two of which are isotopic.
\end{enumerate}

\begin{figure}[htb!]
\def\svgwidth{0.8\textwidth}
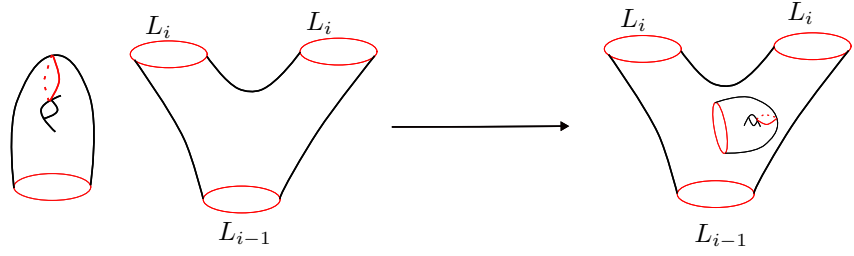
\caption{The pants decomposition $\P_{torus}$ on the left and then the surgery operation on a pair of pants with boundaries at level $i$ and $i-1$.}\label{pantsforpureone}
\end{figure}

We now need to add genus.  Again, we work level by level. Consider $\gamma$ at level $i$ and the induced branch $B_\gamma$ it defines. If $\iota(\E_G)\cap \Ends{B_\gamma}\neq\emp$ we add a punctured torus with a fixed pants decomposition $\P_{\ell(0,1)}$ to the pair of pants or annulus with boundary component $\gamma$ at level $i-1$, see Figure \ref{pantsforpureone}. 

The loops in $\P_{\ell(0,1)}$ get assigned an intermediate level $(i,i+1)$ if they are contained in a region bounded by loops of level $i$ and $i+1$ respectively.

After running this construction for all loops at all levels we obtain a surface $(\Sigma,\P_1)$ with $\P_1$ a collection of separating loops $\Lambda_1$ with a level structure and loops coming from the $\P_{\ell(0,1)}$'s we attached that have an intermediate level structure. It might happen that two loops $\gamma_i,\gamma_j$ in $\Lambda$ at level $i$ and $j$ respectively are isotopic in $S$. If that is the case, we are in an isolated planar end and we delete all but the lowest level loop. This subset of $\Lambda_1$ will be denoted by $\Lambda$ and we will use $\mathcal T$ for the intermediate torus loops. Note that $\P=\Lambda\cup\mathcal T$ still has a natural level structure given by $\Lambda$.

Thus, we obtain a surface $(\Sigma,\P)$ with $\P$ a pants decomposition and $\P=\Lambda\coprod\mathcal T$ and we compactify all isolated planar ends corresponding to the collection $\mathfrak b$, we call the resulting surface $\bar\Sigma$. By construction we have that $\Ends{\bar\Sigma}=\Ends S$, and similarly for the boundary components, and so we have that the surface $\Sigma$ is homeomorphic to $S$, see Theorem \ref{classsurfacethm}.

We now need to show that $\P$ has $(0,1,2)$-biflutes between any two genus ends. Let $x,y\in \E_G\subset \E$. Let $\Sigma_0$ the compact surface cut out by the level one curves. Then $\Sigma_0$ is either a pair of pants, a pair of pants with a punctured torus attached in or a twice punctured torus. Either way $\partial \Sigma_0$ has at least two boundary components that are separating. Thus, there are two distinct components $\gamma_x$, $\gamma_y$ in $\partial\Sigma_0$ such that their branches contain $x$ and $y$ respectively. We then add the surface $P_1(x)$ at the next level that has one boundary component $\gamma_x$ and $P_1(y)$ at the next level that has one boundary component $\gamma_y$. Note that again $P_1(x)$, $P_1(y)$ can be either a pair of pants, a twice or thrice punctured torus.

Continuing this way we get a path of surfaces $\gamma_{xy}\eqdef \cup_{i=1}^\infty P_i(x)\cup P_0\cup_{i=1}^\infty P_i(y)$ such that $P_i(x)$ has one boundary glued to $P_{i-1}(x)$ and $P_1(x)$ has one boundary glued to $P_0$, same for $P_i(y)$.

\textbf{Claim.} Each $P_i(x)$, $P_i(y)$, and $P_0$ is homeomorphic to a pair of pants with a punctured torus glued in.

\begin{figure}[htb!]
\def\svgwidth{0.8\textwidth}
\begingroup%
  \makeatletter%
  \providecommand\color[2][]{%
    \errmessage{(Inkscape) Color is used for the text in Inkscape, but the package 'color.sty' is not loaded}%
    \renewcommand\color[2][]{}%
  }%
  \providecommand\transparent[1]{%
    \errmessage{(Inkscape) Transparency is used (non-zero) for the text in Inkscape, but the package 'transparent.sty' is not loaded}%
    \renewcommand\transparent[1]{}%
  }%
  \providecommand\rotatebox[2]{#2}%
  \newcommand*\fsize{\dimexpr\f@size pt\relax}%
  \newcommand*\lineheight[1]{\fontsize{\fsize}{#1\fsize}\selectfont}%
  \ifx\svgwidth\undefined%
    \setlength{\unitlength}{96.4274453bp}%
    \ifx\svgscale\undefined%
      \relax%
    \else%
      \setlength{\unitlength}{\unitlength * \real{\svgscale}}%
    \fi%
  \else%
    \setlength{\unitlength}{\svgwidth}%
  \fi%
  \global\let\svgwidth\undefined%
  \global\let\svgscale\undefined%
  \makeatother%
  \begin{picture}(1,0.39958893)%
    \lineheight{1}%
    \setlength\tabcolsep{0pt}%
    \put(0,0){\includegraphics[width=\unitlength,page=1]{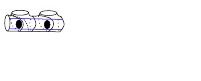}}%
    \put(-0.00043135,0.27747953){\color[rgb]{0.01960784,0,0}\makebox(0,0)[lt]{\lineheight{1.25}\smash{\begin{tabular}[t]{l}$x$\end{tabular}}}}%
    \put(0.32941212,0.27878626){\color[rgb]{0.01960784,0,0}\makebox(0,0)[lt]{\lineheight{1.25}\smash{\begin{tabular}[t]{l}$y$\end{tabular}}}}%
    \put(0,0){\includegraphics[width=\unitlength,page=2]{images/prop56.pdf}}%
    \put(0.45062099,0.38427081){\color[rgb]{0.01960784,0,0}\makebox(0,0)[lt]{\lineheight{1.25}\smash{\begin{tabular}[t]{l}$x$\end{tabular}}}}%
    \put(0.77835677,0.38136198){\color[rgb]{0.01960784,0,0}\makebox(0,0)[lt]{\lineheight{1.25}\smash{\begin{tabular}[t]{l}$y$\end{tabular}}}}%
    \put(0,0){\includegraphics[width=\unitlength,page=3]{images/prop56.pdf}}%
  \end{picture}%
\endgroup%

\caption{How we construct the shift between the ends $x$ and $y$, drawn are a 2-block and a 3-block. The support of the shift is bounded by blue arcs.}\label{prop56claim}

\end{figure}

\bpfc By construction we added a punctured torus to the pair of pants starting at level $i$ if the branch cut off from one the level $i+1$ boundaries contained genus. Since the path is going to the genus ends $x$ and $y$ in both ways this is always the case, see Figure \ref{prop56claim}.

\epfc

Then, by Definition \ref{ABCbiflutedefin} we get that $\gamma_{xy}$ is a $(0,1,2)$-biflute concluding the proof, see Figure \ref{prop56claim}. \epf

We can now prove the main result of the section.

\begin{customthm}{A}\label{puremap} Let $S$ be an infinite-type surface. Then, there exists infinitely many components of Teichm\"uller space, such that $\PMCG S\subset \oMod X$.
\end{customthm} 
\bpf Let $\P$ be the pants decomposition of $S$ coming from Proposition \ref{pantsforpure}. By \cite{PV2018} we know that $\PMCG S$ is the closure of compactly supported mapping classes and handle-shifts. Let $X$ be the hyperbolic structure on $S$ obtained by glueing pants with cuff-lengths equal to one and zero twist. By construction it is a structure of bounded-type.

 Thus, it suffices to show that any handle-shift is in $\oMod X$. Since we have $(0,1,2)$-biflutes between any pair of genus ends this follows from Corollary \ref{bifluteconditionforQC}.
 
To obtain infinitely many components of $\mathcal T(S)$ notice that there are infinitely many $f\in \MCG S$ that are not QC with respect to $X$, see Remark \ref{pAexample}. Let $Y$ be the image structure in $\mathcal T (S)$ and consider $\oMod Y$. Note that, $f\PMCG S f^{-1}\cong \PMCG S$ we can consider maps of the form $fg f^{-1}$ for $g\in \PMCG S$. Since we can realise any $g\in \PMCG S$ as a limit of $g_n\in\oMod X$ we have that $f g_n f^{-1}:Y\rar Y$ is QC and it converges to $fgf^{-1}$ showing that $\PMCG S \subset \oMod Y$. \epf

\brem[Non bounded-type examples] Assume that the surface in Theorem \ref{puremap} contains a genus end $x$ and let $F=F_x$ be either a $(0,1,2)$-biflute or a flute. Let $\gamma_i\subset F$, $i\in\N$, be a sequence of meridians of the handles of $F_i$ that exit one end. Let $X$ be the structure coming from Theorem \ref{puremap} and instead of letting $\ell_X(\gamma_i)=1$ we build a new structure $Y$ for which $\ell_Y(\gamma_i)=\f 1i$ and all other cuffs are unchanged. Then, any shift along $F$ is also modular and so we get components of $\mathcal T(S)$ that are of non bounded-type (as $Y$ has length spectrum accumulating at zero).
\erem
We conclude the Section with a motivating example for Section \ref{REsec}.

\begin{figure}[htb!]
\def\svgwidth{0.8\textwidth}
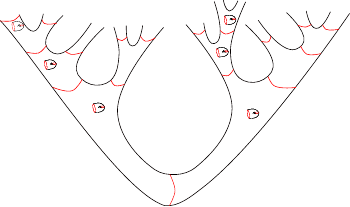
\caption{The Cantor Tree in which we have a genus end labelled by the infinite sequence of left turns $x=\set{LLLLLL\dotsc}$ and an alternating sequence $y=\set{RLRLRLRLRLR\dotsc}$.}\label{figextreegenusone}
\end{figure}

\bese \label{cantortreetwogenus}
Let $S_{x,y}$ be the Cantor Tree surface in which $x\neq y\in\Ends S$ are two genus ends. Assume $S_{x,y}$ looks as Figure \ref{figextreegenusone} and with the given pants decomposition.

This satisfies the condition of Proposition \ref{pantsforpure} and so by Theorem \ref{puremap} the untwisted structure $X$ with bounded pants satisfies  $\PMCG S\subset\oMod X$. Say we want to upgrade the result to any element $\phi\in\MCG S$. Then, essentially, by \cite{dlH2000} we want to generate all tree branch swaps\footnote{A tree branch swap is a (not necessarily level preserving) map that swaps two branches of the tree.}. So we want to cut at two pants curves $\gamma_i$ and $\gamma_{i'}$ and map the branch $B_1=B_{\gamma_i}$ to $B_2=B_{\gamma_{i'}}$ in a QC way.

If $x\in \Ends{ B_1}$ and $y\in \Ends{B_2}$ then it is not obvious how to extend the map mapping pants to pants since $x$ is `straight' and $y$ is alternating. The issue is that we can move each genus in $B_1$ to the appropriate genus in $B_2$ in a QC-way but the supports of these QC-maps are not pairwise disjoint since they all intersect the first pair of pants in $B_1$, therefore is not a priori clear if the limit of these $QC$-maps will even be a homeomorphism.

\begin{figure}[htb!]
\def\svgwidth{0.8\textwidth}
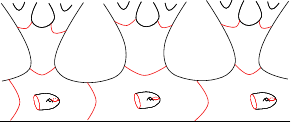
\caption{Another configuration of the Cantor tree with two genus ends obtained by glueing genus to a biflute and Cantor Trees to all boundary components. In this case the genus ends are $x=\set{LLLLLLL\dotsc}$ and $y=\set{RRRRRRR\dotsc}$. }\label{figextreegenustwo}
\end{figure}

However, if we represent $S_{x,y}$ as in Figure \ref{figextreegenustwo}. We have pants curve $\gamma_i(x)$ and $\gamma_i(y)$ cutting out ends neighbourhoods $U_x^i$ and $U_y^i$ that are homeomorphic and have a homeomorphism $\phi_i: U_x^i\rar U_y^i$ mapping pants to pants.

Then, we get that any $\phi\in \MCG{S_{x,y}}$ is in the closure of $\oMod X$, for $X$ the structure obtained glueing pants with cuff-length one and no twist. The homeomorphism from the original Cantor Tree to its ``straightened" version can be realised quasi-conformally on a structure obtained by glueing isometric pants according to the described pants decomposition. This can be seen by a limiting procedure in which at every stage the the Cantor Tree is closer to the ``straight" one, i.e. we swap a genus a time by a tree branch swap.\eese

\brem The condition in Example \ref{cantortreetwogenus} that makes $\oMod X=\MCG S$ will be introduced in Section \ref{Constrsec} and is a condition, that we will call property $\RE$, on homeomorphic ends to have regular neighbourhoods that we can swap via a quasi-conformal homeomorphism of the surface. 
\erem

\section{General Construction}\label{Constrsec}

In this Section we give the definition of property $\RE$ and the main combinatorial and topological constructions needed to show that any surface $S$ with countable end space admits a structure with $\MCG S=\oMod X$.

\subsection{Build Surfaces with Property $\RE$}\label{REsec}
The aim of this section is to carefully define property $\RE$, for regular ends, and build infinite families of infinite-type surfaces with this property. We start with a few definitions.
 
\bdefi\label{REdefi} We say that a pair $(S,\P)$ for $S$ a surface with a pants decomposition $\P$ has \emph{property $\RE$} (regular ends) if:
\begin{enumerate}
\item every two genus ends are connected by a generalised $(A,B,C)$-biflute, if there is a unique genus end $x$ we want a neighbourhood of $x$ to contain a genus $(A,B,C)$-flute;
\item for any two ends $x,y\in \Ends S$ having homeomorphic ends neighbourhoods then, there exists pants curve $\gamma_x,\gamma_y\in\P$ cutting out ends neighbourhoods $U_x$ and $ U_y$ respectively such that:
\begin{itemize}
\item[(2a)] $\phi_{xy}:U_x\rar U_y$ is a homeomorphism; 
\item[(2b)] for $\P_x\eqdef U_x\cap \P$ and $\P_y\eqdef U_y\cap \P$ then $\phi_{xy}(\P_x)=\phi_{xy}(\P_y)$,
\end{itemize}
\end{enumerate}
such neighbourhoods are called $\RE$-neighbourhoods. Moreover, if $S$ has a boundary component $\gamma$ we want any genus ends $x$ to have an $(A,B,C)$-generalised flute containing $x$ as the flute end and with boundary $\gamma$. Finally, if $S$ has a unique genus end $x$ we want an $(A,B,C)$-generalised flute containing $x$.
\edefi 

\bdefi\label{swappableend} Given $(S,\P)$ we say two ends $x$ and $y$ are  \emph{pants-swappable ends} if they have compatible end neighbourhoods as in Definition \ref{REdefi}. The set of swappable ends will be denoted by $\mathcal{SW}(S)=\mathcal {SW}(S,\P)$, where we use the first notation when the pants decomposition is clear.
\edefi

The reason we introduce property $\RE$ is that the map $\phi_{xy}$ maps pants to pants and so, as long its restriction to the pants is uniformly QC and all the pants are glued without twists, it can be realised as a quasi-conformal map, see Lemma \ref{pantsQC}.

\bdefi Given a collection of surfaces $\set{S_i}_{i\in\N}$, with pants decomposition $\P_i$, we call the surface $S=S(\set{S_i}_{i\in\N})$ a \emph{flute of $\set{S_i}_{i\in\N}$ surfaces} if $S$ is obtained by taking a flute $F$ with pants decomposition $\P_F$ and glueing to each $\alpha_i\in\partial F$ the surface $S_i$. Then, $S$ has a pants decomposition $\P=\cup_{i\in\N}\P_i\cup \P_F$. Moreover, we replace each pants in $F$ with a genus pant if one of the boundary components cuts off a genus end.\edefi  

\bese If $S_i$ are punctures then the flute of $S_i$ is just a regular flute. If the $S_i$'s are all punctures except for $S_k$ being a Loch Ness monster then we have that all pants, in $F$, up to level $k$ are genus pants.
\eese 
\bese\label{REex} The second pants decomposition in Example \ref{cantortreetwogenus}, see Figure \ref{figextreegenustwo}, has property $\RE$. Similarly all surfaces $S_i$ of Example \ref{ese3} and flute of flutes of surfaces $S_i$ have property $\RE$, with their induced pants decompositions. Moreover, for all those examples the pants decomposition can be chosen so that the homeomorphisms $\phi_{xy}$ is the ``identity" i.e. they preserve a marking.
\eese

We now state some very useful Lemmas which allows us to build an $\RE$ surface by a flute construction and will allow for transfinite induction arguments.

\blem \label{fluteRE} Let $(S,\P_S)$ be a surface with property $\RE$ and countable end space, a level structure $\Lambda(S)$, and one boundary component which is $\Lambda_0(S)$. Let $S'$ be a flute of $S$ surfaces with $\P_{S'}$ the induced pants decomposition. Then, $(S',\P_{S'})$ has property $\RE$ and a level structure $\Lambda(S)$.
\elem
\bpf 
We label the copies of the surface $S$ that we attached to the flute by $S_i$, $i\in\N$. The end space of $S'$ has a special point $x_\infty$ which corresponds to the flute end and is accumulated by the end spaces of the $S_i$ surfaces. By our assumption on the countable end space and by the fact that all $S_i$ are homeomorphic we have that $x_\infty$ is the unique point of maximal rank, i.e. no other point in $y \in \Ends {S'}$ is such that there exists a neighbourhood $U_y$ homeomorphic to a neighbourhood of $x_\infty$. Thus, we only need to check Definition \ref{REdefi} for points $x,y\in\Ends{S'}\setminus\set{x_\infty}$.

 Then, we have two cases:
\begin{enumerate}
\item $x,y\in\Ends{S_i}$;
\item $x\in\Ends{S_i}$ and $y\in\Ends{S_j}$ with $i\neq j$.
\end{enumerate} 

Case (1) follows from the fact that $S_i$ has property $\RE$, thus we only need to show (2). The surfaces $S_i,S_j$ are equal to $S$ and have the same pants decomposition and we have homeomorphic end neighbourhoods $U_x\cong U_y$. 

If $x$ and $y$ are the same end in $S$ we can take the same pants curve cutting an end neighbourhood $V$. If not we have two points $x,y\in\Ends S$ with homeomorphic ends neighbourhoods so by property $\RE$ we have two pants curve $\gamma_x,\gamma_y$ cutting homeomorphic end neighbourhoods with a homeomorphism sending pants to pants. Pulling back these neighbourhoods to $S_i$ and $S_j$ we obtain the desired neighbourhoods so that $S'$ has property $\RE$.

To obtain the level structure $\Lambda(S)$ we work as follow. Set-wise, we set $\Lambda(S)\eqdef \Lambda(F)\cup \coprod_{i=1}^\infty \Lambda(S_i)$. To define the levels we just translates the levels of the $\Lambda_j(S_i)$. We set:
\[ \Lambda_k(S)=\Lambda_k(F)\cup \cup_{i<k}\Lambda_{k-i}(S_i).\]
If $S$ has genus ends we also need to check the $(A,B,C)$-biflute condition. Let $(A,B,C)$ be the coefficients for the generalised biflute in $S$ and $x\in S_i$ and $y\in S_j$ be genus ends in two distinct copies of $S$, $i<j$. By definition each copy of $S$ has an $(A,B,C)$-flute to the boundary and we can connect them, in $F$ by a chain of genus pants of length $j-i$ by constructions of $F$. This completes the proof. \epf 

If one wants to run similar arguments as before but by glueing non-homeomorphic surfaces $S_i$ one needs a compatibility conditions which we now introduce.

\bdefi \label{compends}Let $(S_1,\P_1)$ and $(S_2,\P_2)$ be surfaces with property $\RE$. We say $(S_1,\P_1)$ and $(S_2,\P_2)$ are \emph{compatible} if $\forall (x,y)\in\Ends{S_1}\times\Ends{S_2}$ with homeomorphic end neighbourhoods there exists:
\begin{itemize}
\item pants curve $\gamma_x\subset \P_1$ and $\gamma_y\subset \P_2$;
\item $\gamma_x$ and $\gamma_y$ cut off ends neighbourhoods $U_x$ and $U_y$;
\item a homeomorphism: $\phi_{xy}:(U_x,U_x\cap \P_1)\rar (U_y,U_y\cap \P_2)$.
\end{itemize}\edefi
Definiton \ref{compends} is just so that $x$ and $y$ have compatible ends neighbourhoods that can be swapped in a quasi-conformal way. Then, for compatible surfaces we obtain a version of Lemma \ref{fluteRE} for connected sums. We first define the induced pants decomposition $\P_{S_1\sharp S_2}$ and level structure $\Lambda(S_1 \sharp S_2)$. We define it as follows. 

\brem[Connected sums, pants decompositions, and levels]\label{ConnSumRe} Let $(S_1,\P_1)$ and $(S_2,\P_2)$ be compatible surfaces with property $\RE$. We first assume they both have at least two genus ends $x_1,y_1$ in $S_1$ and $x_2,y_2$ in $S_2$ (the more complicated case) so we can describe how to glue preserving the $(A,B,C)$-biflute condition. By the $\RE$ condition $S_1$ has an $(A_1,B_1,C_1)$-biflute $F_1$ between $x_1$ and $y_1$. We pick $P_1\in \P_1\cap F_1$ to be a pair of pants with one boundary attached to a 1-handle. We do the same in $S_2$ and get a pants $P_2$. We then glue $S_1$ to $S_2$ by deleting, and identifying the boundaries, of disks in $P_1$ and $P_2$. We then let $\P=\P_1\cup\P_2\cup \gamma \cup\gamma_1\cup\gamma_2$ to be the pants decomposition on $S_1\sharp S_2$ where $\gamma$ is the gluing loop, and $\gamma_1$, $\gamma_2$ are loops in $P_1$, $P_2$ respectively going around a boundary component and $D_1$, $D_2$ respectively. If they only have one genus end we do the same but with the genus flute emanating from $x$.

For the level structure we define set-wise 
\[\Lambda(S_1 \sharp S_2)  \coloneq \Lambda(S_1) \cup \Lambda(S_2) \cup \gamma \cup \gamma_1 \cup \gamma_2\]
 and we'll define $\Lambda_i(S_1 \sharp S_2)$ inductively. We set $\Lambda_0(S_1 \sharp S_2) = \gamma$. Note that this choice is arbitrary, you can choose any level curve to be level $0$. A curve $\lambda \in \Lambda(S_1 \sharp S_2)$ is in $\Lambda_i(S_1 \sharp S_2)$ if there is path connecting $\lambda$ to a curve in $\Lambda_{i-1}(S_1 \sharp S_2)$ that is disjoint from all other level curves.
\erem

Now that we know how to build the pants decomposition for connected sums we can show that property $\RE$ is preserved under connected sums.

\blem\label{sumRElem} Let $(S_1,\P_1)$ and $(S_2,\P_2)$ be compatible surfaces with property $\RE$ then, the connected sum $(S_1\sharp S_2,\P)$, with $\P$ induced by $\P_1$ and $\P_2$, has property $\RE$. Moreover, if $\Lambda(S_1)$ and $\Lambda(S_2)$ are level structures we obtain a new level structure $\Lambda$ on $S_1\sharp S_2$.
\elem
\bpf  We now want to show that $(S_1\sharp S_2,\P)$ has property $\RE$. If two ends are in $S_i$ then it follows from the fact that $S_i$ has property $\RE$. If they are in distinct surfaces then it follows by our compatibility assumption.
 
It remains to check the $(A,B,C)$-biflute condition. If $(S_1\sharp S_2,\P)$ has a single genus end or every genus end is either in $S_1$ or $S_2$ there is nothing to do. So assume we have two genus ends such that $z_1\in S_1$ and $z_2\in S_2$. We now show that they have an $(A',B',C')$-biflute connecting them. 

Let $F_z^1$ be the $(A_1,B_1,C_1)$-biflute connecting $z_1$ to $x_1$ in $S_1$ and let $F_1$ be the $(A_1,B_1,C_1)$-biflute we used to define the connected sum. Since $x_1$ is not an end of $\overline{S_1\setminus F_1}$ we must have that $F_z^1\cap F_1$ is not empty. Then, since they are both biflutes we must have that they overlap. In particular there is a pant $Q$ with a genus handle that is common to both and is the first such pant as seen from $z_1$. Then, we can use $F_1$ to get an $(A_1,B_1,C_1)$-flute from $z_1$ to $\gamma\in\P$. We can do the same in $S_2$ to obtain an  $(A_2,B_2,C_2)$-flute from $z_2$ to $\gamma$. Their, union yields an $(A',B',C')$-biflute connecting $z_1$ to $z_2$ where $A'=\max \set{A_1,A_2}$ and so forth.

The level structure is described in Remark \ref{ConnSumRe}. This completes the proof.\epf

An immediate consequence, which generalizes Lemma \ref{fluteRE}, is as follows.
\bcor \label{sumRE} If we have $n\in\N$ compatible surfaces $(S_i,\P_i)$ then, the connected sum of the $S_i$'s $S$ has a pants decomposition $\P$ such that $(S,\P)$ has property $\RE$ and a level structure $\Lambda$. Similarly, for a flute of compatible surfaces $(S_i,\P_i)$, $i\in\N$ in which there are no ends homeomorphic to the flute end and the $(A_i,B_i,C_i)$ coefficients for the $S_i$'s are uniformly bounded.
\ecor

We now need a topological Lemma to deal with the genus setting of flute surfaces, by a \emph{genus flute} we mean a flute with either finite or infinite genus.

\blem \label{flutehomeo}Let $\set{S_i}_{i\in\N}$ be a collection of surfaces with one boundary component and let $\psi_1,\psi_2:\N\rar\N$ be a bijection and let $F$ be a flute (or genus flute) surface with boundaries labelled by the integers $\geq -1$. Let $\Sigma_1$ and $\Sigma_2$ be the surfaces obtained by glueing to the $i$-th boundary component of $F$ the surface $S_{\psi_j(i)}$, $i\in\N$ and $j=1,2$. Then, there exists a homeomorphism $\phi: \Sigma_1\rar \Sigma_2$ mapping $S_{\psi_1^{-1}(i)}$ to $S_{\psi_2^{-1}(i)}$ via the identity.
\elem 
\bpf Let $\psi\eqdef \psi_2^{-1}\circ\psi_1$ and $f\in\MCG F$ be the map realising the permutation $\psi$ and extend it by the identity to get $\phi:\Sigma_1\rar \Sigma_2$. This is well defined since $\psi$ maps the $i$-th boundary component, glued in $F_1$ to $S_{\psi_1(i)}$ to the $\psi(i)=\psi_2^{-1}\circ\psi_1(i)$ which is glued to $\psi_2(\psi_1(i))$.
\epf

Note that in the above Lemma it is important that we are doing this in a flute and not for example in a ladder surface.

Using Lemma \ref{fluteRE}, Lemma \ref{sumRElem}, and Corollay \ref{sumRE} we will now show that many surfaces enjoy property $\RE$.

\subsubsection{Planar surfaces with countable ordinal type}\label{planarsectionRE}
In this section we show that any planar surface $S$ with countable end space has property $\RE$. This will follow from the fact that closed countable sets of the Cantor sets are given by countable ordinals. Our starting blocks for the transfinite argument will be the $i$-spheres (or spheres with finitely many puncturers), $i\in\N$, of Definition \ref{genflute}.We can now show that any planar surface has property $\RE$.

\bthm\label{planarareRE} Each planar surface $S$ with $\Ends S=\w_0^\zeta n+1\in\cord$ admits a pants decomposition $\P$ with property $\RE$. Moreover, the pants decomposition $\P$ has a sub-collection $\Lambda\subset\P$ of separating loops admitting a level structure.
\ethm 
\bpf By Theorem \ref{homeotypecord} and Corollary \ref{sumRE} it suffices to show that each surface $S_{\zeta,n}$, $\Ends{S_{\zeta,n}}\cong \w_0^\zeta n+1$, has property $\RE$. Moreover, since $S_{\zeta,n}$ is the connected sum of $n$ copies of $S_{\zeta}$ by Lemma \ref{sumRElem} it suffices to prove it for $S=S_\zeta$ i.e. a $\zeta$-sphere, and then take connected copies with the same pants decomposition so that they satisfy the compatibility condition of Lemma \ref{sumRElem}.

Let $\w\in\cord$ be simple, i.e. $\w=\w_0^\zeta+1$. Then, we need to build $S_\w$, that is the surface with $\Ends S=\w$. By Lemma \ref{simplecord} we can do this recursively by starting with finite ordinals. Each finite ordinal is a surface of type $S_n=S_{(0,n)}$ where $(0,n)$ is the Cantor-Bendixon rank of $\Ends{S_n}$ and they are the surfaces of Example \ref{ese3} which have property $\RE$ by Example \ref{REex}. Then, by Lemma \ref{fluteRE} the first countable ordinal $\w_0$ has property $\RE$. Thus, by Corollary \ref{sumRE}, since at every stage we have compatible surfaces, and transfinite induction we get that all $S_{\w}$ have property $\RE$.

Moreover, the level structure comes from the fact that the $S_n$ surfaces have a natural level structure and the flute as well. Then if we glue a surface $S_n$ to level, $i$ in the flute, we just translate the $S_n$ levels by $i$. This, can be carried on by the transfinite induction. \epf

\brem Theorem \ref{planarareRE} with Theorem \ref{mainthm} already gives us uncountably many examples of surfaces in which $\oMod X\cong\MCG S$. We will also give examples with uncountable end space in Section \ref{uncendsexample}.
\erem

\subsection{Countable surface with genus ends}\label{Regenussec}

We now show that every surface $S$ with countable end space has property $\RE$, i.e. the analogue to Theorem \ref{planarareRE} but we now allow genus ends. Since we already dealt with planar surfaces and if $g(S)<\infty$ we can connect sum a finite genus surface to get the required statement, we only need to deal with surfaces with $g(S)=\infty$.

Let $(\E, \E_g)$ be a countable end space, as before the aim will be to build a surface $S=S(\E, \E_g)$ with:
\begin{itemize}
\item a pants decomposition $\P$;
\item a level structure $\Lambda$;
\end{itemize}
such that $(S,\P)$ has property $\RE$.

In this setting, in contrast to Corollary \ref{planarremarkRE}, we do not have obvious models for every end type then the idea is to build a system of ends neighbourhoods $\mathcal U\eqdef \set{U_{\alpha,\beta},V_{\alpha,\beta}}_{\alpha,\sim\beta\in \E}$ with $\alpha\in U_{\alpha,\beta}$, $\beta\in V_{\alpha,\beta}$, $U_{\alpha,\beta}\cong V_{\alpha,\beta}$, and satisfying certain compatibility conditions (for example any two elements in $\mathcal U$ are either nested or disjoint). We then coherently endow the $U_{\alpha,\beta}$'s with compatible pants decompositions to guarantee property $\RE$. Since $\mathcal U$ covers all ends we have that $\overline {S\setminus\mathcal U}$ is a compact surface and so we can extend the pants decomposition on $\mathcal U$ to one on the whole of $S$.

Let $S$ be a model for our surface $S=S(\E,\E_g)$. Since our surface has countably many ends we can work using the derived sets, see Definition \ref{cantorbendixonrankdefi}. We will denote by $D_g^\alpha=D^\alpha( \E)\cap \E_g$, for $\alpha$ an ordinal, the genus ends of order $\alpha$. Then, the $x,y\in D_g^\alpha$ that have homeomorphic ends neighbourhoods in $S$, $x\sim y$, partition $D_g^\alpha$ ($D^\alpha$) into equivalence classes that we denote by $\hat D_g^\alpha=\set{[x]}$ ($\hat D^\alpha$) where each class has at most countably many elements and there are at most countably many classes.

We start with a Corollary which follows from the previous section.

\bcor\label{planarremarkRE} Let $S$ be a surface with countable end space and let $x$ be a planar end of rank $\alpha$. Then, $x$ has a closed neighbourhood homeomorphic to the $\alpha$-sphere $S_\alpha$ with an open disk removed.
\ecor

 This is the main definition of the Section and it will allow us to construct our $\RE$ neighbourhoods.
\bdefi\label{normalneighs} Given a surface $S=\cup_{k\in\N} S_k$ with end space $\E\cong \w_0^\zeta n +1$, we say that a system of closed end neighbourhoods $\U\eqdef \U^p\cup\U^g$:
\[ \U^p\eqdef \set{U^p_{x,y},V^p_{x,y}}_{x\sim y\in\E\setminus \E_g}, \U^g\eqdef\set{U_{x,y}^g,V^g_{x,y}}_{x\sim y\in \E_g}, \]
is \emph{normal} if it satisfies:
\begin{enumerate}
\item each $U\in\U^p$ is planar and each $U\in\U^g$ is infinite-genus and they both have a single boundary component;
\item $x\sim y\in \mathcal D^\alpha$,$*=p,g$:
\begin{itemize} 
\item[(2a)] $\rank{U_{x,y}^*}=\rank x=\alpha=\rank y=\rank{V_{x,y}^*}$ and $x$, $y$ are the unique such points;
\item[(2b)] $x\in U_{x,y}^*$ and $y\in V_{x,y}^*$ are pairwise disjoint, i.e. $U^*_{x,y}\cap V^*_{x,y}=\emp$;
\item[(2c)] there is a homeomorphism $\phi_{xy}:U_{x,y}^*\rar V_{x,y}^*$;
\item[(2d)] for $[x]\in\hat D^\alpha$ and $[x]=\set{x,y_i}_{i\in\N}$ then $U_{x,y_j}^*\subset U_{x,y_i}^*$ for all $j\geq i$;
\end{itemize}
\item for all $\alpha<\beta\leq\zeta$ and $x\sim y \in \mathcal D^\alpha$, $z\sim w\in\mathcal D^\beta$,$*=p,g$ we have:
\begin{itemize}
\item[(3a)] $U_{x,y}^*,V_{x,y}^*$ and $U_{z,w}^*$, $V_{z,w}^*$ are either disjoint or the $\alpha$-neighbourhoods are contained in the $\beta$ ones;
\item[(3b)] $U_{x,y}^*,V_{x,y}^*$ are contained in at most finitely many $U_{z,w}^*$, $V_{z,w}^*$;
\end{itemize}
\item for every $k\in\N$ there are finitely many $U_{x,z}^*,V_{x,y}^*$ intersecting $S_k$;
\item let $x\in\E_g$ be a genus end and consider the neighbourhood $U_{x,y}^*$ for $y\in [x]$, then: 
\begin{itemize} 
\item[(5a)] there exists a maximal open covering $W$ of the ends of $U_{x,y}^g$ distinct from $x$ such that $U'_x=\overline {U_{x,y}^g\setminus W}$ has infinite genus, same for the $V_{x,y}$ sets;
\item[(5b)] if $U_{x,y}^g\in \U^g$ so is $\hat U_{x,y}^g\cong  U_{x,y}^g$ contained in $ U_{x,y}^g$ with $\overline{ U_{x,y}^g\setminus \hat  U_{x,y}^g}$ a twice punctured torus, same for the $V_{x,y}$ sets.
\end{itemize}
\end{enumerate} 
\edefi  
\brem Condition (2b) fails once we allow perfect sets in the end space, for example in the case of a Cantor Tree.\erem

In the next Lemma, for $\xi\in\cord$, we will denote the collection of end neighbourhoods $\U_{<\xi}$ given by $\U^p\cup \U^g_{<\xi}$, where in $\U^g_{<\xi}$, we only have neighbourhoods for genus ends of rank $<\xi$. This will be needed for the version of Lemma \ref{planarnormalneighs} containing genus, in particular to deal with property (5a) of Definition \ref{normalneighs}.

\blem\label{maximalcovering} Let $S$ be a surface of infinite-type, $\E\cong \w_0^\zeta n+1$, with an exhaustion $\set{S_k}_{k\in\N}$, and a system of end neighbourhoods $\U_{<\xi}$, $\xi\leq \zeta$. Then, for $U$ a neighbourhood of a genus end $x$ of $\rank x=\rank U =\xi$ there is a neighbourhood $W\subset\U_{<\xi}$ of $\Ends U \setminus\set x$ such that $W$ is not contained in any other collection of elements of $\U_{<\xi}$.
\elem 
\bpf We will use Zorn's Lemma \cite{C1997} in which if $A,B$ are end neighbourhoods covering $\Ends U \setminus\set x$ then $A\leq B$ if $A\subset B$. Thus, we need to show that chains admit maximal elements. 

Let $\mathcal C\eqdef \set{ C_i}_{i\in I}$ be a chain of end neighbourhoods with $C_i \leq C_{i+1}$. Let $A\in\pi_0( C_i)=\set{A_i^j}_{j\in\N}$ be a component and let $S_k$ be large enough so that $\partial A\subset S_k$ is an essential loop\footnote{That means that $\partial A$ is not homotopic into a power of a boundary component or bounds a disk.}. Then, if $A_j\subset C_j$, $j\geq i$, is the component containing $A$ it must intersect $S_k$. However, by property (4) of Definition \ref{normalneighs}, there are finitely many such elements in $\U_{<\xi}$ and so we can extract a maximal neighbourhood containing $A$. 

Specifically we start by considering $A_1^1$. Then, by the above there exists $i_1\in \N$ such that, up to relabelling:
\[ \forall i\geq i_1: A_1^1\subset A_i^1,\qquad A_1^j=A_1^k, k\geq i_1.\]
We then, repeat this for $A_{i_1}^2$ and so forth. Let $C_M\eqdef \cup_{j\in\N} A^j_{i_j}$, we claim it is maximal, i.e. $\forall i :~C_i\leq C_M$. This, follows because for all $i$ $\pi_0( C_i)=\set{A_i^j}_{j\in\N}$ and for each $j$ there exists $i_j$ such that $A^j_i\subset A^j_{i_j}\subset C_M$.\epf

We now show that a system of closed end neighbourhoods $\U$ always exists for planar surfaces, that is we build $\U^p$ in Definition \ref{normalneighs}. For $\alpha\in\cord$ we will use $\U^p_{<\alpha}$ to denote a system of end neighbourhoods as in Definition \ref{normalneighs} but only having neighbourhoods for points of ranks $<\alpha$.

\blem\label{planarnormalneighs} Given a surface $S=S(\E,\E_g)$ with $\E\cong \w_0^\zeta n+1$ and a compact nested exhaustion $\set{S_k}_{k\in\N}$ of $S$, then there exists a normal system of planar closed end neighbourhoods $\mathcal U^p$.
\elem 
\bpf By Corollary \ref{planarremarkRE} we have that any planar end $x\in D^\alpha$ has a standard neighbourhood obtained by glueing surfaces $\Sigma_i$ to a flute $F$ where $\rank{\Sigma_i}=\alpha_i<\alpha$ and $\alpha_i\rar\alpha$, in the case of successor ordinal the $\alpha_i$ are all the same and for limit ordinals they are distinct and $\alpha_i<\alpha_{i+1}$. Note that we can ignore property (5) of Definition \ref{normalneighs} as we are dealing with planar ends. Let $\eta\in\cord$, $\eta\leq \zeta$, be the smallest ordinal such that any planar end has rank at most $\eta$. The proof of the Lemma will be by transfinite induction.

 We define $\U^p_0$ to be the collection of pairwise disjoint end neighbourhoods corresponding to isolated planar ends such that each $U\in \U^p_0$ is a punctured disk. These clearly satisfy (1)-(4) since they are always pairwise disjoint and properly embedded in $S$. Since these are all the isolated points this is the base case of our transfinite induction and we now deal with the induction step.

Assume that we built neighbourhoods satisfying (1)-(4) for all planar ends with rank less than $\alpha$ that is, we have defined $ \U^p_{<\alpha}$. We now prove the transfinite induction step.

\textbf{Claim.} Given $ \U^p_{<\alpha}$ we can build $ \U^p_{\leq\alpha}$.

\bpfc Pick $x\in D^\alpha(\E\setminus \E_g)$ then, by Corollary \ref{planarremarkRE} there is a unique equivalence class, i.e. we have that $\forall x,y\in D^\alpha(\E\setminus \E_g)$ then $x\sim y$. We label the ends as $\set{x,y_i}_{i\geq 1}=D^\alpha(\E\setminus \E_g)$. 

Pick a planar neighbourhood $U_x$ such that $\partial U_x$ has a single component that is disjoint from all elements of $\U^p_{<\alpha}$, and such that $x$ is the unique point of $\rank x =\alpha$. That is $U_x$ is a $\alpha$-sphere by Corollary \ref{planarremarkRE} and $\overline{U_x\setminus \U_{<\alpha}^p}$ is a flute with boundary. Moreover, since every end in $U_x$ has rank strictly less than $\alpha$ we have a similar neighbourhood $V_{xy_i}$, see Corollary \ref{planarremarkRE}. These neighbourhoods can be built by using Lemma \ref{maximalcovering}. By taking the same neighbourhoods $U_{xy_i}=U_x$ for all $y_i$ we can extend $ \U^p_{<\alpha}$ to $ \U^p_{\leq\alpha}$. The neighbourhoods $U_{xy_i}$ and $V_{xy_i}$ so defined clearly satisfy (1), (2), and (3). Property (4) again follows by the maximal rank neighbourhoods being pairwise disjoint and the $S_k$ being compact.
\epfc

By the claim we can apply transfinite induction to obtain the required result.\epf

We now do the transfinite induction step for the genus version of Lemma \ref{planarnormalneighs}, which will build us $\U^g$ and so $\U$.

\blem\label{systemofcompneighsstep} Given $S=S(\E,\E_g)$ with $\E\cong \w_0^\zeta n+1$ and a compact nested exhaustion $\set{S_k}_{k\in\N}$ of $S$ and a system of normal neighbourhoods $\U_{<\xi}$, $\xi\leq \zeta$, then we can enlarge $\U_{<\xi}$ to $\U_{\leq\xi}$.
\elem 
\bpf The proof will be a series of inductions. The most important ones being showing that such a collection of neighbourhoods exists for a single class $[x]\in \hat D^\xi( \E_g)$. Then, we repeat it for each class.

 Let $[x]\in \hat D^\xi( \E_g)=\set{[x_j]}_{j\in\N}$, and label each end in the class as $[x]=\set{x,y_i}_{i\geq 1}$. For ease of notation we will use $U_i$, $V_i$ instead of $U_{x,y_i}$, $V_{x,y_i}$ and first prove the Lemma for a single class.

By definition of $x\sim y_1$, we have neighbourhoods $U_1$, $V_1$ and a homeomorphism $\phi_1: U_1\rar V_1$. Note that (1) is automatically satisfied since we are dealing with genus ends and up to restricting the neighbourhoods we can assume that they have a single boundary component and that are pairwise disjoint (2b-c). By using the fact that the end space is countable, again by restricting the neighbourhoods if necessary, we can assume that $U_1$, $V_1$ also satisfy (2a). Thus, we can assume that $U_1$ and $V_1$ satisfy (1) and (2a-c) of Definition \ref{normalneighs}.

We now modify $U_1$, $V_1$ to obtain (3a-b) while preserving (1) and (2a-c). This follows by the fact that if we restrict neighbourhoods we do not break properties (1), (2a-c) of Definition \ref{normalneighs}.

\textbf{Property (3a) and (3b).} The simple closed curves $\partial U_1$ and $\partial V_1$, intersect at most finitely many elements of $\U_{<\xi}$ by property (4). We denote this collection by $H$. Let 
\[E_1\eqdef \left(H\cup \phi_1^{-1}\left( H\cap V_1\right)\right)\cap U_1,\]
and note that $\overline {U_1\setminus E_1}\subset U_1$ contains a neighbourhood $U_1'$ of $x$ and $U_1'$ is homeomorphic to the $y_1$ neighbourhood $\phi_1\left(U_1'\right)\subset V_1$. With an abuse of notation we still denote these new neighbourhoods by $U_1$ and $V_1$. Thus, $U_1$ and $V_1$ now satisfy (1), (2a-c), and (3a). Moreover, (3b) follows by induction.

\vspace{0.5cm}

We now deal with properties (2d), (4), and (5) of Definition \ref{normalneighs}.

\vspace{0.5cm}

\textbf{Property (5a).} By Lemma \ref{maximalcovering} we have a maximal covering $W_x\subset U_1$ of ends distinct from $x$. Then, $\hat U_1=\overline{U_1\setminus W_x}$ is an infinite-type surface of either finite or infinite genus. If $\hat U_1$ is of infinite genus there is nothing to do. Otherwise we have infinitely many elements of $W_x$ of infinite-genus, using property (5a) we have a covering $W'_x\subset W_x$ such that $\overline{W_x\setminus W'_x}$ is a countable disjoint union of twice punctured tori. Thus, $\hat U_1=\overline{U_1\setminus W'_x}$ has infinite-genus. The same can be done for $V_1$ with a covering $W_y$. Thus, both $\hat U_1$ and $\hat V_1$ are homeomorphic to infinite genus flutes with boundary.

\textbf{Property (5b).} Since each neighbourhood is of infinite-genus this follows by classification.

\textbf{Property (2d) and (4).} To verify these properties we first need to construct all other neighbourhoods $U_i$, $V_i$, $i>1$. Since there are countably many of those indexed by $\N$ we can work by induction, with base case $i=1$. 

Assume we defined $U_j$, $V_j$ up to $i-1$ satisfying properties (1)-(5). We now build $U_i$ and $V_i$. Since $x\sim y_i$ we have $U_i\cong V_i$ and up to restricting them we can assume that $U_i\subset U_{i-1}$, giving (2d). Moreover, we can also assume them to be pairwise disjoint and, as before, to satisfy (1) and (2a-c). We then repeat the proof of (3a), (5a), (5b) to conclude. Property (4) follows from the nestedness condition.

This concludes the proof for a single class in $\hat D^\xi(\E_g)$.

\vspace{0.5cm}

The rest of the argument can be done by usual induction as the classes of $\hat D^\xi(\E_g)$ are indexed by $\N$. The base case is the proof we just did. We define $\U_{\xi,n}$ to be the collection of closed ends neighbourhoods containing $\U_{\xi}$ and the neighbourhoods covering the ends in the classes up to $[x^n]$. Then, to go from $n$-classes to $n+1$-classes we start with $U_{\xi,n}$ and we want to cover $[x^{n+1}]$. Then, we repeat the above steps but now using $\U_{\xi,n}$ instead of $\U_\xi$.\epf

The proof of Proposition \ref{systemofcompneighs}, which constructs the neighbourhoods $\U$, will be by transfinite induction. We first show the base case.

\blem\label{basecasesystemofcompneighs} Given $S=S(\E,\E_g)$ with $\E\cong \w_0^\zeta n+1$ and a compact nested exhaustion $\set{S_k}_{k\in\N}$ of $S$. Then, there exists a system of normal ends neighbourhoods $\mathcal U_{<1}$.\elem 
\bpf We begin by taking $\mathcal U$ to be the collection of planar end neighbourhoods $\mathcal U^p$ of Lemma \ref{systemofcompneighs} and add to it pairwise disjoint neighbourhoods of all isolated genus ends homeomorphic to a Loch Ness monster minus a disk. Then, these satisfy (1)-(5) of Definition \ref{normalneighs}. We now need to deal with genus ends that are isolated in $\E_g$ but are accumulated by planar ends. 

Consider $\hat D_0^g=\set{[x^k]}_{k\in\N}$ with each class given by $[x^k]=\set{x^k, y_i^k}_{i\in\N}$. Each element has a neighbourhood that is homeomorphic to a genus flute $F_x$  with planar surfaces glued in, the homeomorphism between them is the one of Lemma \ref{flutehomeo} applied to a maximal covering of $\Ends F \setminus\set x$, as in Lemma \ref{maximalcovering}, by elements of $\U^p$. By picking them to be pairwise disjoint and so that their boundaries avoid elements of $\U^p$ we are done. \epf

We can finally show that surfaces with countable end space admit a normal system of ends neighbourhoods as in Definition \ref{normalneighs}.

\bprop\label{systemofcompneighs}  Given $S=S(\E,\E_g)$ with $\E\cong \w_0^\zeta n+1$ and a compact nested exhaustion $\set{S_k}_{k\in\N}$ of $S$. Then, there exists a system of normal ends neighbourhoods $\mathcal U\eqdef \U^p\cup\U^g=\set{U_{x,y},V_{x,y}}_{x\sim y\in\E}$.
\eprop 
\bpf  We induct on the genus rank. By Lemma \ref{basecasesystemofcompneighs} we have $\U_{<1}$. By transfinite induction and Lemma \ref{systemofcompneighsstep} we obtain the required system $\U=\U_{\leq\zeta}$.\epf

We now have to build the pants decomposition $\P$ on $S$ on which to test property $\RE$, see Definition \ref{REdefi}. The neighbourhoods will be given by the elements of $\U$ and we only need to guarantee the $(A,B,C)$-flute condition.

We will do this by taking compatible pants decomposition in the elements of the system of normal ends neighbourhood $\U$. We will start with lower order ones and then take nested simple multi-curves to obtain a pants decomposition on $\cup_{U\in\U}U$.

The end space $\E$ of $S$ can be decomposed as follows:
\[ \E=\E_g\cup \E_1\coprod\E_2,\] 
where the planar ends $\E_2$ do not accumulate on genus ends. Then, by Theorem \ref{planarareRE} we know that all ends in $\E_1\cup\E_2$ have $\RE$ neighbourhoods. Let $\gamma_0$ be a simple closed loop in $S$ splitting $\E_2$ from $\E_g\cup\E_1$ and let $\mathcal U^p\subset\mathcal U$ be the subset of neighbourhoods of the planar ends. Let $\P^p$ be the simple multi-curve obtained by adding to $\gamma_0$ the pants decomposition $\P^p_0$ of Theorem \ref{planarareRE} in each set of $\mathcal U^p$ starting from smaller rank to bigger rank. Then, in $(S,\P^p)$ we have that each pair of planar ends $x\sim y$ have neighbourhoods $U_x$ and $U_y$ and a homeomorphism $\phi_{xy}:\left(U_x,\P^p\vert_{U_x}\right)\rar\left(U_y,\P^p\vert_{U_y}\right)$, i.e. $\RE$ neighbourhoods.

\bdefi\label{multicurvedefi}
Given a system of end neighbourhoods $\U$ on $S$ a simple multi-curve $\P^\alpha$ inducing pants decomposition on elements of $\U_{<\alpha}$ is said to be an\emph{ $\RE$ multi-curve} if the elements of $\U_{<\alpha}$ are $\RE$-neighbourhoods with $(0,1,2)$ generalised flutes, see Definition \ref{REdefi}. We say $\P$ is an \emph{$\RE$-system} if it is an $\RE$ multi-curve on $\U=\U_{\leq \zeta}$, $\Ends S\cong \w_0^\zeta n+1$.

\edefi 

We now deal with genus ends and work inductively on the rank, where the base case is given by isolated genus ends. The next Lemma will be the induction step of the argument and Lemma \ref{planarREsystem} will be the base case.

\blem\label{inductgenusRE} Given $S$ with a system of end neighbourhoods $\U$ and an $\RE$ multi-curve $\P^\alpha$ on $\U_{<\alpha}$ then, we have a multi-curve $\P^{\bar\alpha}$ containing $\P^\alpha$ that is an $\RE$ multi-curve on $\U_{\leq \alpha}$.
\elem 
\bpf Consider $U_{x,y}$ and $V_{x,y}$, $x\sim y$ of rank $\alpha$, in a system of normal neighbourhoods $\U$. By property (3) of Definition \ref{normalneighs} we can add their boundaries to $\P^\alpha$. 

Let $\E_x$ and $\E_y$ be the ends of $U_{x,y}$ and $V_{x,y}$ respectively and $\E_x'\eqdef \E_x\setminus\set{x}$ and $\E_y'\eqdef \E_y\setminus\set{y}$. Then, by property (5a) and hypothesis $\E'_x$ and $\E'_y$ admit maximals covers $W_x$, $W_y$ by pairwise-disjoint $\RE$-neighbourhoods. Then, we define $U'_{x,y}\eqdef \overline{U_{x,y}\setminus W_x}$ and $V_{x,y}'\eqdef \overline{V_{x,y}\setminus W_y}$. Moreover, note that $\partial U_{x,y}',\partial V_{x,y}'\subset \P^\alpha$.

By (5a) of Definition \ref{normalneighs} we get that the $U_{x,y}'$ and $V_{x,y}'$ are both either infinite genus or planar flutes. Then, we extend the multi-curve $\P^\alpha$ on $U_{x,y}$ to $U_{x,y}'$ so that $U'_{x,y}$ is a $(0,1,2)$-flute. Since by hypothesis each genus end in $W_x$ has a $(0,1,2)$-flute to $\partial W_x$ it can be glued to obtain a $(0,1,2)$-generalised flute to $\partial U_{x,y}$. Then, to get $\P^\alpha_x$ we push it to $V_{x,y}$ via the homeomorphism of Lemma \ref{flutehomeo}. 

Repeating the above construction for every equivalence class in $D_g^\alpha$, choosing a representative in each one, and using the fact that the $U_{x,y}$ and $V_{x,y}$ are pairwise disjoint yields $\P^{\bar\alpha}$. \epf

We first build a multi-curve that forms an $\RE$ system for planar ends.

\blem\label{planarREsystem}  Let $S$ be a surface with countable end space. Then, there exists a simple multi-curve $\P^p$ on $S$ and a system $\U^p$ of closed end neighbourhoods, covering all planar ends, forming an $\RE$-system.
\elem 
\bpf Let $\xi$ be the maximal rank of a planar end. The argument will be by transfinite induction. We start first by taking $\P^0$ to be the simple multi-curve cutting out all isolated planar ends. This clearly satisfies the condition in Definition \ref{multicurvedefi}. We then induct on the rank of the planar ends and use Lemma \ref{inductgenusRE} to conclude by taking $\P^p\eqdef \P^\xi$.
\epf

By another transfinite induction argument we can extend the result to genus ends.

\bcor\label{REsystem} Let $S$ be a surface with countable end space. Then, there exists a simple multi-curve $\P$ on $S$ and a system of closed end neighbourhoods $\U$ forming an $\RE$-system.
\ecor 
\bpf We take the closed end neighbourhoods $\U=\U^p\cup \U^g$ of Proposition \ref{systemofcompneighs} and now build the multi-curve $\P$. The construction is by transfinite induction on the rank with inductive step Lemma \ref{inductgenusRE}.

By Lemma \ref{planarREsystem} we only need to deal with genus ends since we have a planar $\RE$-system $\P^p$. For each genus end $x$, isolated in $\E^g$, we have a neighbourhood $U_x\in\U^g$ such that $x\in U_x$ is the unique genus end any other end is planar and thus covered by $\U^p$. Then, the non-compact component $F_x$ of $\overline{U_x\setminus \P^p}$ is a genus flute with boundary and end given by $x$. We then, endow $F_x$ with the pants decomposition $\P_x$ turning it into a $(0,1,2)$-flute. Then, for any other end $y\in [x]\in \hat D_g^0$ we push the multi-curve $\P_x$ with the homeomorphism of Lemma \ref{flutehomeo}. This yields a simple multi-curve $\P_x$ extending $\P^p$ so that if forms an $\RE$-system for any planar end and any end in $[x] \in\hat D_g^0$. Repeating this for all $x\in D_g^0$ completes the base case and the transfinite induction step is Lemma \ref{inductgenusRE}. \epf

Thus, we finally have:
\bthm\label{cendsareRE} Any surface $S$ with countable end space has property $\RE$. Moreover, the pants decomposition $\P$ has a level structure $\Lambda$.
\ethm 
\bpf By Corollary \ref{REsystem} we have an $\RE$-system with simple multi-curve $\P$. Since $\mathcal U$ is a covering of the end space we have that $\Sigma_0\eqdef \overline{S\setminus \U}$ is compact. Then, we can extend the simple multi-curve $\P$ to a pants decomposition $\P$ which can be given a level structure by choosing a root in $\Sigma_0$.  
\epf

\subsection{For $\RE$ surfaces we have $X\in\mathcal T(S)$ with $\MCG S=\oMod X$.} \label{REMCGsec}

We now conclude this work by showing that surfaces with property $\RE$ satisfy $\oMod X=\MCG S$ for infinitely many structures $X\in\mathcal T(S)$ in different components of $\mathcal T(S)$. The last Lemma is the following.

\blem\label{endspaceapp} Let $(S,\P)$ be a surface with property $\RE$, $X$ an untwisted bounded structure on $S$, $\phi\in\MCG S$, and $\Sigma\subset S$ a compact sub-surface with no finite-genus complementary regions. Then, there exists a quasi-conformal homeomorphism $f:X\rar X$ such that:
\[f\vert_\Sigma=\phi\vert_\Sigma\]

\elem
\bpf 
The Lemma will follow once we show the following claim.

\textbf{Claim 1.} There exists $x_1,\dotsc, x_n\in\mathcal E$ and neighbourhoods $U_i\subset\Sigma^C$ of $x_i$ and a homeomorphism $\psi$ such that:
\begin{itemize}
\item $\partial U_i\subset\P$, $U_i\cap U_j=\emp$ and $\psi(\partial U_i)\subset \P$;
\item for $y_i\eqdef \phi(x_i)$ and $V_i\eqdef \psi(U_i)\subset\phi(\Sigma)^C$ we have a homeomorphism $\psi_i: (U_i,\P\cap U_i)\rar (V_i,V_i\cap \P)$ and $\psi: (U,\P\cap U)\rar (V,V\cap \P)$ for $U\eqdef \cup_{i=1} U_i$ and $V\eqdef \cup_{i=1}^n V_i$.
\end{itemize} 
\bpfc By property $\RE$ for every $x\in\mathcal E$ there are neighbourhoods $U_x$, $V_{\phi(x)}$ of $x$ and $\phi(x)$ respectively such that $\partial U_x,\partial V_{\phi(x)}$ are in $\P$ and a homeomorphism $\psi_x: (U_x,U_x\cap \P)\rar (V_{\phi(x)},V_{\phi(x)}\cap\P)$. Moreover, since $\Sigma$ is compact we can assume that each $U_x\subset \Sigma^C$ and $V_\phi(x)\subset \phi(\Sigma)^C$. Then, the collection of open sets $\set{U_x\times V_{\phi(x)}}_{x\in\mathcal E}$ is an open covering of the compact set $\mathcal E$ so it has a minimal finite covering of the form $U_{x_i}\times V_{\phi(x_i)}$ for $1\leq i\leq n$. We let $y_i\eqdef \phi(x_i)$, $U_i=U_{x_i}$, $V_i\eqdef V_{y_i}$.

Then, since the $U_i$'s cover $\mathcal E$ and $\phi$ is a bijection we must have that the $V_i$ cover $\E$. Since $\partial U_i\subset\P$ we have that $\partial U_i\cap \partial U_j=\emp$ and so $U_i$ and $U_j$ are either disjoint or nested. Since we took a minimal finite cover they must be disjoint. Moreover, still by property $\RE$, the $U_i$ and $V_i$ come equipped with homeomorphisms:
\[\psi_i:(U_i,U_i\cap \P)\rar (V_i,V_i\cap\P),\]
then we define $\psi: \coprod_{i=1}^n(U_i,U_i\cap \P)\rar\coprod_{i=1}^n (V_i,V_i\cap\P)$ by $\psi\vert_{U_i}\eqdef \psi_i$.\epfc

Let $U\eqdef \cup_{i=1}^n U$, $V\eqdef \cup_{i=1}^n V_i$ and $\psi: (U, U\cap \P)\rar (V,V\cap \P)$ be as above. Then, since $X$ is untwisted and $\psi$ maps pants to pants we have that $\psi$ is a quasi-conformal map, see Lemma \ref{pantsQC}. We now need to extend $\psi $ to the whole surface. 

\begin{figure}[htb!]
\def\svgwidth{0.8\textwidth}
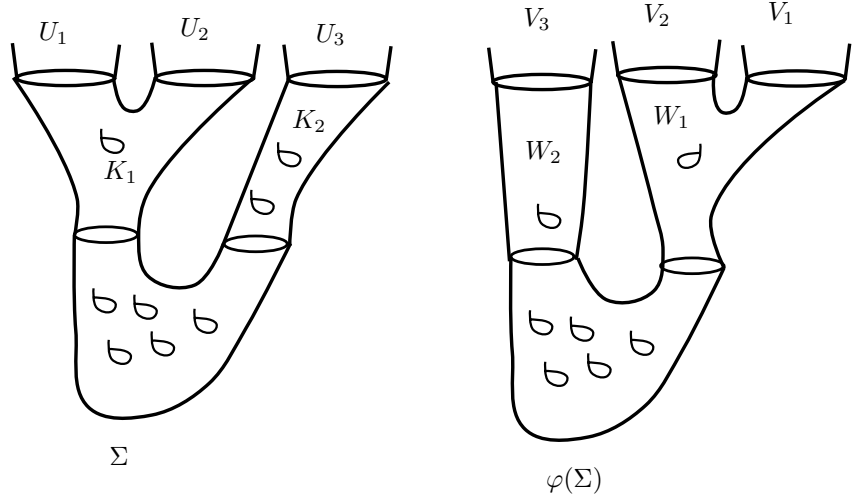
\caption{The two decomposition of $S$ into $K=K_1\cup\Sigma\cup K_2$, $U_1$, $U_2$, and $U_3$ and the one with $W$'s and $V$'s. }\label{figtwodecomps}
\end{figure}

Let $K\eqdef \overline{S\setminus U}$ and $W\eqdef \overline{S\setminus V}$ and notice that by construction $\Sigma\subset K$ and $\phi(\Sigma)\subset W$. Also by construction if $U_i$ is contained in $H\in\pi_0( \Sigma^C)$ then $V_i$ is contained in $\phi(H)\in\pi_0( \phi(\Sigma)^C)$, this is because $H$ is a neighbourhood of $U_i$ and so needs to be mapped to a neighbourhood of $V_i$. Thus, complementary components of $\overline{K\setminus \Sigma}=\set{K_j}$ and $\overline{W\setminus\phi(\Sigma)}=\set{W_j}$ can be paired up by their boundaries and differ by at most genus, see Figure \ref{figtwodecomps}. Moreover, if a component $K_j$ is planar and does not cut out genus ends then $W_j$ is also planar, and vice-versa since otherwise there would be a complementary component of $\Sigma $ of finite genus.

In particular, $K_j$ and $W_j$ need to have the same number of boundary components. Thus, if $W_j$ and $K_j$ are not homeomorphic they differ by genus and cut out at least a genus end, see Figure \ref{figtwodecomps}.

\textbf{Claim 2.} Given, $U_i$, or $V_i$, as above with infinite-genus and $k\in\N$ there exists a quasi-conformal embedding $s_i^k:U_i\hookrightarrow U_i$ such that $\overline{U_i\setminus \eta_i^k(U_i)}\cong S_{k,2}$ (or $t_i^k$ for $V_i$).

\bpfc Since $U_i$ has infinite-genus by property $\RE$ we have that there exists an $(A,B,C)$-flute $F$ connecting $\partial U_i$ to a genus end. Consider the shear $s_i: U_i\rar U_i$ along $F$ as in Lemma  \ref{boundedfluteshiftismodular}. This is a quasi-conformal map such $\overline{U_i\setminus s_i(U_i)}$ is a twice punctured torus. Then, the $k$-th composition $s^k$ of $s$ is quasi-conformal embedding $s_i^k:U_i\hookrightarrow U_i$ with $\overline{U_i\setminus s_i^k(U_i)}\cong S_{k,2}$.\epfc

Then, by doing finitely many shears we can make the $W_j$ and $K_j$ homeomorphic and so we get a homeomorphism $f_0:K\rar W$ extending $\phi$. On each of the ends neighbourhoods we define the map by 

\[\tilde\psi_i\eqdef t_i^{h}\circ\psi_i\circ\left(s_i^k\right)^{-1}:s_i^k(U_i)\rar V_i,\]

which is a QC-homeomorphism, note that $h,k$ can be zero in which case the maps are the identity. We then define $f_1=\cup _{i=1}^n\tilde \psi_i$ and so the glueing of $f_0$ an $f_1$ along $\partial K$ yields the desired quasi-conformal homeomorphism $f:X\rar X$ with $f\vert_{\Sigma}=\phi\vert_{\Sigma}$.\epf

Given Lemma \ref{endspaceapp} we can finally prove our main result.

\begin{customthm}{C}\label{mainthm} Let $(S,\P)$ have property $\RE$ then there exists infinitely many component of $\mathcal T(S)$ such that for all $X\in \mathcal T(S)$ we have $\MCG S=\oMod X$.
\end{customthm}
\bpf First note that it suffices to show it for a hyperbolic structure $X$ and this automatically gives the statement for the whole Teichm\"uller component  $\mathcal T(X)$ containing $X$.

Given $(S,\P)$ let $X$ be the structure obtained by glueing pants with cuff lengths equal to $1$ and no twists. Let $\set{S_n}_{n\in\N}$ be the compact exhaustion of $S$ obtained by cutting $S$ at level $n$ and taking $S_i$ to be the unique connected compact component. By enlarging, if needed, the compact elements of the exhaustion we can also assume that for all $n$ no complementary region of $S_n$ has finite genus.

Let $\phi\in\MCG S$, we now want to show that we can find maps $\phi_n\in \oMod X$ such that $\phi_n\rar \phi$ in the compact open topology. Since each $S_n$ is so that no complementary region has finite genus by Lemma \ref{endspaceapp} we have quasi-conformal homeomorphisms $\phi_n$ agreeing with $\phi$ on $S_n$ and so we conclude the first part of the Theorem.
To get infinitely many components of the Teichm\"uller space we proceed as in Theorem \ref{puremap}. This concludes the proof. \epf

By combining Theorem \ref{mainthm} and Theorem \ref{cendsareRE} we obtain our remaining main result.

\begin{customthm}{B} Given any infinite-type surface $S$ with countably many ends there exists a hyperbolic structure $X$ such that $\MCG S = \oMod X$.
\end{customthm}

\subsection{Surface with non-countable end spaces}\label{uncendsexample}

We conclude this work by showing that property $\RE$ is not restricted to surfaces with countable end space. The examples we now build will all be based on trees for simplicity.

\bdefi
We define a \emph{standard tree of type $\w\in\cord$} as a Cantor Tree with the tree pants decomposition such that in each pants we attach a planar surface of type $\w$. We denote it by $\triangle(\w)$. Similarly we define the genus version $\triangle_g(\w)$ in which we start with a blooming Cantor-Tree.
\edefi

\brem\label{stdtreeareRE} Note that for $\w\neq \w'$ then any trees $\triangle(\w)$, $\triangle(\w')$, $\triangle_g(\w)$, $\triangle _g(\w')$ satisfy, if equipped with the natural pants decomposition, the conditions of Lemma \ref{sumRE} and so have property $\RE$.
\erem 

Thus, we obtain that:
\blem Let $\w_1,\dotsc,\w_n,\w_1',\dotsc,\w_m'\in\cord$ be distinct ordinals and $S'$ be any surface with countable end space. Let $S$ be the surface obtained by doing a connected sum of the $\triangle(\w_i)$, $\triangle_g(\w_j')$, and $S'$. Then, $S=S(\vec\w,\vec \w',S')$ has property $\RE$. 
 \elem 
\bpf By Theorem \ref{cendsareRE} and Remark \ref{stdtreeareRE} all these surfaces have property $\RE$. Because $\w_i\neq\w_j$ and $\w_i'\neq\w_j'$ we get that he only homeomorphic ends of $\triangle(\w_i)$, $\triangle_g(\w_j')$, and $S'$ are the planar ends. Since these are all canonically homeomorphic and so we can assume them to have compatible ends with respect to their pants decompositions. Then, by Corollary \ref{sumRE} we have the stated result. \epf 

\brem The surfaces $S(\vec\w,\vec \w',S')$, where $S'$ is a surface with countable end space, $\vec\w\in\cord^n$, $\vec w'\in\cord^m$, form another uncountable family. Moreover, as long as one picks genus ends in a compatible way, see Corollary \ref{sumRE}, we can also glue genus ends of countable rank to the tree to obtain an equivalent result.\erem 
 
  \brem In this Remark we show that the surface of \cite{MR2022} has property $\RE$. Specifically the surfaces built for \cite[Theorem 1.1]{MR2022} which states \emph{there exist examples of surfaces that have non self-similar end spaces with a unique maximal type of end and set of maximal ends homeomorphic to a singleton or to a Cantor set.}

Mann and Rafi do two constructions, one for uncountable planar end space and one for countable genus ends. The genus case follows by Theorem \ref{mainB} and the examples in the planar uncountable case can also be given an $\RE$ structure with the techniques developed in this work.
 \erem

 \section{Questions}
Based on our examples we wonder the following natural questions.

\textbf{Question 1.} For every surface $S$ there exists a pants decomposition $\P$ such that $S$ has property $\RE$?
 
 Note that in all our arguments of the previous section we intensively used the structure of countable sets of $\Cantor$ and being able of finding maximal rank elements. This is not true for general subsets and the countable rank points form an open set accumulating on a perfect set $K\cong\Cantor$. Thus, finding standard neighbourhoods of points $x\in K$ might require analysing all possible ways that finite rank points can accumulate onto $x$. 
 
 \textbf{Question 2.} Can one use modular maps to obtain some kind of Nielsen-Thurston classification?
 
 \nocite{C20171,C2018c,CS2018}

\thispagestyle{empty}
{\small
\markboth{References}{References}

\bibliographystyle{alpha}
\bibliography{mybib}{}

}
	\end{document}